\makeatletter

\documentclass[12pt,a4paper]{article}
\usepackage[fleqn]{amsmath}
\usepackage{amscd,amssymb}
\usepackage{amsthm}
{
\theoremstyle{plain}
\newtheorem{th@@r@m}{Theorem}[section]
\newtheorem{l@mm@}[th@@r@m]{Lemma}
\newtheorem{pr@p@s@t@@n}[th@@r@m]{Proposition}
\newtheorem{c@r@ll@ry}[th@@r@m]{Corollary}
}
{
\theoremstyle{definition}
\newtheorem{d@f@n@t@@n}[th@@r@m]{Definition}
\newtheorem{r@m@rk}[th@@r@m]{Remark}
\newtheorem{@x@mpl@}[th@@r@m]{Example}
\newtheorem{pr@bl@m}[th@@r@m]{Problem}
}
\newif\ifInsection

\renewcommand{\section}{\@startsection {section}{1}{\z@}%
                                   {-28pt \@plus -3pt \@minus -1pt}%
                                   {16pt \@plus1pt}%
                                   {\normalfont\large\bfseries}}

\newdimen\setsub@nh@b@
\def\@seccntformat#1{%
\hbox to \setsub@nh@b@{\S\csname the#1\endcsname.}\ignorespaces
}

{
\theoremstyle{plain}
\newtheorem*{MTh@@r@m}{Main Theorem}
\newtheorem*{TH@@r@m}{Theorem}
}

{
\theoremstyle{definition}
\newtheorem*{REM@rk}{Remark}
}

\newenvironment{thm}{\Inthmtrue\begin{th@@r@m}}{\end{th@@r@m}\Inthmfalse}
\newenvironment{lem}{\Inthmtrue\begin{l@mm@}}{\end{l@mm@}\Inthmfalse}
\newenvironment{prop}{\Inthmtrue\begin{pr@p@s@t@@n}}{\end{pr@p@s@t@@n}\Inthmfalse}
\newenvironment{cor}{\Inthmtrue\begin{c@r@ll@ry}}{\end{c@r@ll@ry}\Inthmfalse}
\newenvironment{df}{\Inthmtrue\begin{d@f@n@t@@n}}{\end{d@f@n@t@@n}\Inthmfalse}
\newenvironment{rem}{\Inthmtrue\begin{r@m@rk}}{\end{r@m@rk}\Inthmfalse}
\newenvironment{ex}{\Inthmtrue\begin{@x@mpl@}}{\end{@x@mpl@}\Inthmfalse}

\newenvironment{pf}{\Inthmtrue\begin{proof}}{\end{proof}\Inthmfalse}

\newenvironment{claim}[1][Claim]{
   \trivlist
  \thm@preskip\topsep
  \thm@postskip\thm@preskip
  \@topsep \thm@preskip               
  \@topsepadd \thm@postskip           
   \item[\hskip\labelsep\itshape#1\@addpunct{.}]\ignorespaces
}{\endtrivlist}

\newenvironment{pfofclaim}[1][Proof of Claim]{\par
  \normalfont
  \topsep6\p@\@plus6\p@ \trivlist
  \item[\hskip\labelsep\itshape
    #1\@addpunct{.}]\ignorespaces
}{%
  \endtrivlist
}

\newcounter{subthm}[th@@r@m]
\newif\ifInthm
\def\MultiNumbering{
 \def\Inthmlabel{\theth@@r@m.\thesubthm}
 \def\c@equation{\ifInthm\c@subthm\else\c@th@@r@m\fi}
 \def\theequation{\ifInthm\Inthmlabel\else\theth@@r@m\fi}
}

\def\TagsOnLeft{\tagsleft@true}
\def\TagsOnRight{\tagsleft@false}

\def\joken{\refstepcounter{subthm}%
  \edef\@currentlabel{\theth@@r@m.\thesubthm}%
     \ignorespaces{\@currentlabel}}

\evensidemargin=\oddsidemargin
\def\SkipAmount{
\baselineskip=17pt plus 2pt
\lineskiplimit=3pt
\lineskip=3pt plus 1pt
}

\def\rangai{%
	\newbox\@currbox
	\newbox\@marbox
  \ifhmode
    \@bsphack
    \@floatpenalty -\@Mii
  \else
    \@floatpenalty-\@Miii
  \fi
  \ifinner
    \@parmoderr
    \@floatpenalty\z@
  \else
    \@next\@currbox\@freelist{}{}%
    \@next\@marbox\@freelist{\global\count\@marbox\m@ne}%
       {\@floatpenalty\z@
        \@fltovf\def\@currbox{\@tempboxa}\def\@marbox{\@tempboxa}}%
  \fi
	\@rangaiypar
}
\long\def\@rangaiypar#1{%
  \@savemarbox\@marbox{#1}%
  \global\setbox\@currbox\copy\@marbox
  \@rangaixypar
}
\long\def\@savemarbox #1#2{%
  \global\setbox #1%
    \color@vbox
      \vtop{%
        \hsize\marginparwidth
  			\parindent\z@ \parskip\z@skip
 			 \linewidth\hsize
 			 \@totalleftmargin\z@
 			 \leftskip\z@skip \rightskip\z@skip \@rightskip\z@skip
			  \parfillskip\@flushglue \lineskip\normallineskip
			  \baselineskip\normalbaselineskip
			  \sloppy		
        \reset@font				         
        \normalsize
        #2%
        \outer@nobreak
        }%
    \color@endbox
}
\def \@rangaixypar{%
  \ifnum\@floatpenalty <\z@\@cons\@currlist\@marbox\fi
  \setbox\@tempboxa
    \color@vbox
      \vbox \bgroup							
      \par\vskip\z@skip      
      \outer@nobreak
    	\egroup                  
  	\color@endbox																																
  	\ifnum\@floatpenalty <\z@
    	\@largefloatcheck
    	\@cons\@currlist\@currbox
    	\ifnum\@floatpenalty <-\@Mii
     	 \penalty -\@Miv
      	\@tempdima\prevdepth
     	 \vbox{}%
      	\prevdepth\@tempdima
      	\penalty\@floatpenalty
    	\else
      	\vadjust{\penalty -\@Miv \vbox{}\penalty\@floatpenalty}\@Esphack
    	\fi
  	\fi
  \@esphack
}

\newcommand{\Twoside}{\@twosidetrue  \@mparswitchtrue}

\let\bib\bibitem              
\def\lb#1{\label{#1}}

\newif\ifInFtn
\newif\ifFtnMark
  
  \def\NoFootnoteMark{\FtnMarkfalse}

\renewcommand{\@makefntext}[1]{%
  \parindent 1em%
  \ifFtnMark
  	\noindent\hb@xt@1.8em{\hss\@makefnmark}#1
   \else
   	\noindent\hb@xt@1.8em{}#1
    \fi}
 
\def\ShowTagFlags{
  \newbox\hid@rikakko       \newbox\m@gikakko
  \newdimen\haba@temp     \newdimen\migi@temp
  \setbox\hid@rikakko\hbox{(}     \setbox\m@gikakko\hbox{)}
  \haba@temp=\textwidth
  \advance\haba@temp\marginparsep
  \migi@temp=\marginparsep
  \advance\migi@temp\wd\m@gikakko
  
  \def\lbp##1{\label{##1}\rule{0pt}{0pt}\marginpar{\texttt{##1}}\ignorespaces}
  
	\def\lbs##1{\label{##1}\rangai{\texttt{\normalsize##1}}\@afterheading}

	\def\footl@b@l{}
  \def\lbfn##1{\label{##1}\gdef\footl@b@l{\texttt{##1}}}
 
  \long\def\@footnotetext##1{\insert\footins{%
      \reset@font\footnotesize
      \interlinepenalty\interfootnotelinepenalty
      \splittopskip\footnotesep
      \splitmaxdepth \dp\strutbox \floatingpenalty \@MM
      \hsize\columnwidth \@parboxrestore
      \protected@edef\@currentlabel{%
         \csname p@footnote\endcsname\@thefnmark
      }%
      \InFtntrue
      \color@begingroup
        \@makefntext{%
          \rule\z@\footnotesep\ignorespaces##1\hfill\rlap{\hskip\marginparsep\footl@b@l}\@finalstrut\strutbox}
      \color@endgroup}
      \let\footl@b@l\@empty\InFtnfalse}
  
 \def\lbe##1{\label{##1}}

 \def\tagform@##1{%
 	\iftagsleft@%
  	 \maketag@@@{\rlap{\hskip\haba@temp\texttt{\df@label}}%
 																														(\ignorespaces##1\unskip\@@italiccorr)}%
   \else%
   	\maketag@@@{(\ignorespaces##1\ignorespaces
    	 \rlap{\hskip\migi@temp\texttt{\df@label}}\unskip\@@italiccorr)}%
   \fi
 }

  \let\lb\relax
  \def\lb##1{\ifmmode\lbe{##1}
  										\else
            						\ifInFtn
                  				\lbfn{##1}
                      	\else
                       	\ifInthm
                       		\lbp{##1}
                         \else
                         	\ifInEnum
                         		\lbp{##1}
                           \else
                           	\lbs{##1}
                            \fi
                         \fi
                       \fi
                     \fi}
  
  \def\bib{\@ifnextchar[{\xbib}{\ybib}}
    \def\xbib[##1]##2{\bibitem[##1]{##2}\rule{0pt}{0pt}\marginpar{\texttt{##1}}\ignorespaces}
    \def\ybib##1{\bibitem{##1}\rule{0pt}{0pt}\marginpar{\texttt{##1}}\ignorespaces}
}

\def\citex[#1]#2{[\texttt{#2, #1}]}
\def\citey#1{[\texttt{#1}]}

\def\KillRefCite{
  \def\ref##1{\texttt{##1}}
  \def\eqref##1{\texttt{(##1)}}
  \DeclareRobustCommand\cite{\@ifnextchar[{\citex}{\citey}}
}

\newbox\q@dsymb@l
\setbox\q@dsymb@l\hbox{\qedsymbol}
  \newdimen\q@dsk@p   
  \q@dsk@p=\textwidth
  \advance\q@dsk@p-\wd\q@dsymb@l
\def\qedinmmode{
\iftagsleft@\tag*{\rlap{\hskip\q@dsk@p\qedsymbol}}
\else\tag*{\qedsymbol}
}

\def\killqed{\renewcommand{\qed}{}}

\def\same{{\hbox to 1cm {\hrulefill}} }

\def\cond[#1]{\par\noindent\rlap{(#1)}\hskip35pt
\hangindent=35pt\hangafter=1}

\def\indno#1{\vskip\itemsep\vskip\parsep%
                       \par\hskip1.2em\llap{\upshape(#1)}\hskip.5em\ignorespaces}
\def\no#1{{\hskip.5em\upshape(#1)}\hskip.5em\ignorespaces}

\def\dkajoindent{0pt}
\def\dkajo{%
\vskip1ex\par\noindent
\hangindent=\dkajoindent%
\hangafter=1%
\llap{\hbox to 1em{\hfil$\bullet$\hfil}}%
\hskip\dkajoindent%
}

\def\RHeads#1{\markright{\hfill\textsl{#1}\hspace{2em}}}

\newif\ifInEnum

\def\parenumerate{%
	\ifnum \@enumdepth >\thr@@\@toodeep\else
    \advance\@enumdepth\@ne
      \edef\@enumctr{enum\romannumeral\the\@enumdepth}%
      \expandafter
      \list
        \csname label\@enumctr\endcsname
			  {\setlength{\leftmargin}{0pt}\setlength{\labelwidth}{1.2em}
         \setlength{\labelsep}{.5em}
		     \addtolength{\itemindent}{\labelwidth}
         \addtolength{\itemindent}{\labelsep}
         \addtolength{\itemindent}{\parindent}
         \setlength{\topsep}{\itemsep}
         \setlength{\parsep}{0pt}
         \usecounter\@enumctr\def\makelabel##1{\hss\llap{\upshape##1}}}%
  \fi}

\renewcommand{\enumerate}{%
  \ifnum \@enumdepth >\thr@@\@toodeep\else
    \advance\@enumdepth\@ne
    \edef\@enumctr{enum\romannumeral\the\@enumdepth}%
      \expandafter
      \list
        \csname label\@enumctr\endcsname
			  {\setlength{\leftmargin}{0pt}\setlength{\labelwidth}{1.2em}
         \setlength{\labelsep}{.5em}
		     \addtolength{\leftmargin}{\labelwidth}
         \addtolength{\leftmargin}{\labelsep}
         \addtolength{\leftmargin}{\parindent}
         \setlength{\topsep}{\itemsep}
         \setlength{\parsep}{0pt}
           \usecounter\@enumctr\def\makelabel##1{\hss\llap{##1}}}%
  \fi}

\newenvironment{parenumr}{\InEnumtrue

	\begin{parenumerate}}{\end{parenumerate}\InEnumfalse}

\newenvironment{parenuma}{\InEnumtrue

	\begin{parenumerate}}{\end{parenumerate}\InEnumfalse}

\def\Thm#1{Theorem #1}

\def\Cor#1{Corollary #1}

\def\Prop#1{Proposition #1}

\def\Lem#1{Lemma #1}
\def\Lems#1{Lemmas #1}
\def\Df#1{Definition #1}

\def\Rem#1{Remark #1}

\def\case#1{\par\noindent{\it Case} {#1}.\quad\ignorespaces }
\def\qedclaim{\hfill $\lozenge$}

\def\bnm#1#2{\binom{#1}{#2}}
\def\bnmt#1#2{{\tbinom{#1}{#2}}}

\def\bnmm#1#2{\bnm{#1}{#2}_+}
\def\bnmmt#1#2{\bnmt{#1}{#2}_+}

\def\ch#1{\operatorname{char}(#1)}

\def\ee#1#2#3{#1 \leq #2 \leq #3}
\def\ne#1#2#3{#1 < #2 \leq #3}
\def\en#1#2#3{#1 \leq #2 < #3}

\def\dss#1#2{\mathop\bigoplus_{#1}^{#2}}
\def\sm#1{\mathop\sum_{#1}}      

\def\smm#1#2{\mathop\sum_{#1}^{#2}}

\def\un#1{\mathop\bigcup_{#1}}

\def\ep#1{\mathop{\bigwedge}^{\kern-0.3pt#1}\;}

\def\Hsupb#1#2#3{H^{#1}_{#2}(#3)}

\def\Hm#1#2{H^{#1}_\mathfrak{m}(#2)}

\def\Hn#1#2{H^{#1}_\mathfrak{n}(#2)}

\def\Hom#1#2#3{\text{\rm Hom}_{#1}(#2,#3)}
\def\Homcal#1#2#3{{{\mathcal H}\kern-1.2pt\textit{om}}\,_{#1}(#2,#3)}

\def\Tor#1#2#3#4{\hbox{\rm Tor}_{#1}^{#2}(#3,#4)}  
\def\Ext#1#2#3#4{\operatorname{Ext}^{#1}_{#2}(#3,#4)}  
\def\Extcal#1#2#3#4{{{\mathcal E}\kern-1pt\textit{xt}}\,^{#1}_{#2}(#3,#4)}
\def\Ker#1{\text{\rm Ker}(#1)} 
\def\Kersup#1#2{\text{\rm Ker}^{#1}(#2)} 
\def\Ima#1{\operatorname{Im}(#1)}
\def\Imasup#1#2{\operatorname{Im}^{#1}(#2)}
\def\Coksup#1#2{\operatorname{Coker}^{#1}(#2)}

\def\Syz#1#2#3{\operatorname{Syz}^{#2}_{#1}(#3)}

\def\Spec#1{\operatorname{Spec}(#1)}

\def\GL{\textit{GL}}
\def\Pfaff{\textrm{Pfaff}}
\def\gradplain#1#2{[#1]_{#2}}
\def\gradn#1#2{\lt[#1\rt]_{#2}}
\def\grad#1#2{\bigl[#1\bigr]_{#2}}

\def\ann#1{\operatorname{ann}(#1)}
\def\directlim#1{
\ooalign
{\hfil\hbox{$\varinjlim$}\hfil\crcr
\leavevmode\hfil\raise-10pt\hbox{\scriptsize$#1$}\hfil}
}

\def\len#1#2{l_{#1}(#2)}

\def\depm#1{\operatorname{depth}_\mfrak(#1)}

\def\hei#1{\operatorname{ht}(#1)}
\def\rk#1#2{\operatorname{rank}_{#1}(#2)}
\def\rank{\operatorname{rank}}

\def\dimk#1{\operatorname{dim}_k(#1)}

\def\socd#1{\operatorname{socd}(#1)}
\def\shmin#1{\#_{\textit{min}}(#1)}

\def\row#1#2#3{{({#1}_{#2},{\ldots}\,,{#1}_{#3})}} 

\def\trow#1#2#3{{{}^t\kern-2pt({#1}_{#2},{\ldots}\,,{#1}_{#3})}}
\def\tpose#1{{{}^t\kern-1pt#1}}
\def\trowiv#1#2#3#4{{{}^t\kern-1.2pt({#1},{#2},{#3},{#4})}}

\def\seq#1#2#3{{{#1}_{#2},{\ldots}\,,{#1}_{#3}}} 
\def\seqsub#1#2#3{{{#1}_{#2},{\ldots},{#1}_{#3}}}

\def\ddd{,{\ldots}\,,}

\def\op{\oplus}

\def\ot{\otimes}

\def\tm{\times}    
\def\eqp{\overset{\text{\rm{p}}}{=}}
\def\opdot{{
\,\ooalign
{\hfil\raise 0.35ex\hbox{$\centerdot$}\hfil\crcr
\leavevmode\hbox{$\oplus$}}\,
}}

\def\indexonsmile#1{
\ooalign
{\hspace{-3pt}\raise -0.3ex\hbox{\scriptsize$\smallsmile$}\crcr
\leavevmode\hfil\raise0.45ex\hbox{\scriptsize$#1$}\hfil}
}

\def\hra{\hookrightarrow}
\def\thra{\twoheadrightarrow}
\def\lra{\longrightarrow}
\def\xr@ght@rr@w[#1]#2{\xrightarrow[#1]{#2}}
\def\xxr@ght@rr@w#1{\xrightarrow{#1}}
\def\xra{\@ifnextchar[{\xr@ght@rr@w}{\xxr@ght@rr@w}}
\def\xl@ft@rr@w[#1]#2{\xleftarrow[#1]{#2}}
\def\xxl@ft@rr@w#1{\xleftarrow{#1}}
\def\xla{\@ifnextchar[{\xl@ft@rr@w}{\xxl@ft@rr@w}}

\def\la{\langle}
\def\ra{\rangle}

\def\bt{\bullet}
\def\cc{\circ}
\def\lt{\left}
           
\def\rt{\right}
\def\eset{\emptyset}
\def\sset{\subset}

\def\bslash{\backslash}

\def\Rstar{{R^{\star}}}
\def\Rbarstar{{{\bar R}^{\star}}}

\def\Wsup#1{{W^{#1}}}
\def\Xsub#1{{X_{#1}}}
\def\Xbarsub#1{{{\bar X}_{#1}}}
\def\thgsub#1{{\thg_{#1}}}
\def\Ysub#1{{Y_{#1}}}
\def\Ybarsub#1{{{\bar Y}_{#1}}}
\def\thgsub#1{{\thg_{#1}}}
\def\phgsub#1{{\phg_{#1}}}
\def\Usub#1{{U_{#1}}}

\def\Nbf{\text{\bf N}}

\def\Pbf{\text{\bf P}}

\def\Zbf{\text{\bf Z}}

\def\lbs{\boldsymbol{l}}

\def\mfrak{\mathfrak{m}}
\def\nfrak{\mathfrak{n}}

\def\pfrak{\mathfrak{p}}

\def\Kcal{\mathcal{K}}

\def\Atild{{\tilde A}}

\def\Etild{{\tilde E}}
\def\Ftild{{\tilde F}}

\def\Htild{{\tilde H}}
\def\Itild{{\tilde I}}

\def\Ztild{{\tilde Z}}

\def\btild{{\tilde b}}

\def\etild{{\tilde e}}
\def\ftild{{\tilde f}}
\def\gtild{{\tilde g}}

\def\vtild{{\tilde v}}

\def\hgtild{{\tilde \hg}}

\def\kgtild{{\tilde \kg}}

\def\Thgtild{{\tilde \Thg}}

\def\tAtild{{^t\kern-0.2em \Atild}}
\def\tA{{^t\kern-0.2em A}}
\def\tB{{^t\kern-0.2em B}}
\def\tF{{^t\kern-0.2em F}}
\def\tFtild{{^t\kern-0.2em \Ftild}}
\def\tPhg{{^t\kern-0.1em\Phg}}
\def\Jlsup#1{{{}^{#1}\kern-0.2em J}}
\def\Kcallsup#1{{{}^{#1}\kern-0.1em\Kcal}}

\def\Abar{{\bar A}}

\def\Fbar{{\bar F}}
\def\Gbar{{\bar G}}
\def\Hbar{{\bar H}}
\def\Ibar{{\bar I}}
\def\Jbar{{\bar J}}

\def\Lbar{{\bar L}}
\def\Mbar{{\bar M}}

\def\Rbar{{\bar R}}

\def\Zbar{{\bar Z}}

\def\ebar{{\bar e}}
\def\fbar{{\bar f}}
\def\gbar{{\bar g}}
\def\hbar{{\bar h}}

\def\nbar{{\bar n}}

\def\vbar{{\bar v}}
\def\wbar{{\bar w}}

\def\ngbar{{\bar \ng}}

\def\mfrakbar{{\bar \mfrak}}

\def\lgbar{{\bar \lg}}

\def\ngbar{{\bar \ng}}

\def\pgbar{{\bar \pg}}

\def\fhat{{\hat f}}
\def\ghat{{\hat g}}

\def\Pcheck{{\Check P}}

\def\gcheck{{\Check g}}

\def\pcheck{{\Check p}}

\def\ag{\alpha}
\def\bg{\beta}
\def\dg{\delta}
\def\gg{\gamma}
\def\kg{\kappa}
\def\Lg{\varLambda}
\def\lg{\lambda}
\def\mg{\mu}
\def\ng{\nu}
\def\sg{\sigma}
\def\tg{\tau}
\def\eg{\varepsilon}
\def\ig{\iota}
\def\pg{\pi}
\def\vpg{\varpi}
\def\thg{\theta}

\def\Thg{\varTheta}

\def\Gg{\varGamma}

\def\Phg{\varPhi}
\def\Psg{\varPsi}
\def\og{\omega}
\def\Og{\varOmega}

\def\Xg{\varXi}
\def\xg{\xi}
\def\rg{\rho}
\def\phg{\varphi}
\def\psg{\psi}
\def\hg{\eta}
\def\chg{\chi}

\def\mod{\operatorname{mod}}

\def\MATr{\text{\rm MAT}}

\def\fort{\text{for}}

\def\andt{\text{and}}

\def\andr{\text{\rm and}}

\def\fallt{\text{for all}}

\def\CM{Cohen-Macaulay}

\def\Gor{Gorenstein}

\def\Wei{Weierstrass}

\def\fgg{finitely generated graded}

\def\Grob{Gr\"obner}

\newenvironment{smallmatrixl}{\null\,\vcenter\bgroup
  \Let@\restore@math@cr\default@tag
  \baselineskip8\ex@ \lineskip1.5\ex@ \lineskiplimit\lineskip
  \ialign\bgroup$\m@th\scriptstyle##$\hfill&&\thickspace
  $\m@th\scriptstyle##$\hfill\crcr
}{%
  \crcr\egroup\egroup\,%
}
\newenvironment{matrixl}{%
  \hskip -\arraycolsep\array{*\c@MaxMatrixCols l}%
}{%
  \endarray \hskip -\arraycolsep
}

\newenvironment{smallbmatrix}{\lt[\begin{smallmatrix}}
{\end{smallmatrix}\rt]}

\newcommand{\cdotsfor}[1]{%
  \ifx[#1\@xp\scdots@for\else\cdots@for\@ne{#1}\fi}
\def\scdots@for#1]{\cdots@for{#1}}
\def\cdots@for#1#2{\multicolumn{#2}c%
  {\m@th\dotsspace@1.5mu\mkern-#1\dotsspace@
   \xleaders\hbox{$\m@th\mkern#1\dotsspace@\cdot\mkern#1\dotsspace@$}%
           \hfill
   \mkern-#1\dotsspace@}%
}
  
\def\mat#1#2#3#4{
\lt(\begin{matrix} #1 & #2 \vspace{5pt} \\  
          #3 & #4 
   \end{matrix}\rt)
}

\def\mattsgn#1\emattsgn{
\text{\footnotesize$\lt(\begin{matrix} #1\end{matrix}\rt)$}
}

\def\bclmsc#1#2{
\lt[ \smallmatrix  {#1} \\ 
{#2}
\endsmallmatrix \rt]
}

\def\supbr#1#2{{{#1}^{[#2]}}}
\def\supang#1#2{{{#1}^{\langle#2\rangle}}}
\def\supbrp#1#2{{{#1}^{[#2]}{}'}}
\def\supangp#1#2{{{#1}^{\langle#2\rangle}{}'}}

\def\suplpar#1#2{{{}^{(#2)}\kern-1pt{#1}}}
\def\sublpar#1#2{{{}_{(#2)}\kern-1pt{#1}}}
\def\BR#1{B_R(#1)}
\def\BRp#1{B_{R'}(#1)}

\newcommand{\kernsub}[1]{\kern-0.2em{{}_{#1}}}
\newcommand{\kernnsub}[1]{\kern-0.3em{{}_{#1}}}

\def\Pthree{{\Pbf^3}}

\def\Gbarbt{\Gbar_\bt}

\def\Gbt{G_\bt}

\def\Gsub#1{G_{#1}}
\def\Gpsub#1{G'_{#1}}
\def\Gppsub#1{G''_{#1}}

\def\overcc#1{{\ooalign{\hfil\raise1.7ex\hbox{\scriptsize$\cc$}\hfil\crcr $#1$}}}

\def\overtildcc#1{{\ooalign{\hfil\raise2.2ex\hbox{\scriptsize$\cc$}\hfil\crcr $#1$}}}

\def\lgbt{\lg_\bt}

\def\lgbar{{\bar\lg}}

\def\lgbarbt{{\lgbar_\bt}}

\edef\today{\ifcase\month\or
  January\or February\or March\or April\or May\or June\or
  July\or August\or September\or October\or November\or December\fi
  \space\number\day, \number\year}

\def\@maketitle{%
  \newpage
  \null
  \vskip 2em%
  \begin{center}%
    {\LARGE \@title \par}%
    \vskip 1.5em%
    {\large
      \lineskip .5em%
      \begin{tabular}[t]{c}%
        \@author
      \end{tabular}\par}%
    \vskip .5em%
    {\footnotesize \@date}%
  \end{center}%
  \par
  \vskip 1.5em}

\def\Watashidesu{Mutsumi A{\scshape masaki}\vspace{1ex}\\\tokoro\\\email}

\renewcommand{\alignedat}[2][c]{%
    \start@aligned{#1}{#2}
}

\makeatother

\def\same{{$\rule{15mm}{0.4pt}$\,}}

\def\Asix{{\it Application of the generalized Weierstrass preparation 
theorem to the study of homogeneous ideals}, 
Trans. AMS {\bf 317} (1990), 1 -- 43.}

\def\Aeleven{{\it Generators of graded modules associated with
linear filter-regular sequences},
J. Pure Appl. Algebra {\bf 114} (1996), 1 -- 23.}

\def\Aseventeen{{\it Generic \Grob\ bases and \Wei\ bases of homogeneous
submodules of graded free modules}, 
J. Pure Appl. Algebra {\bf 152} (2000), 3 -- 16.}

\def\Atwenty{{\it Inequalities satisfied by the basic sequence
of an integral curve in $\Pthree$}, arXiv.org e-Print archive, math.AC/0504131,
April 2005.}

\def\HUl{{\it General hyperplane sections of algebraic varieties},
J. Alg. Geom. {\bf 2} (1993), 487 -- 505.}

\def\HMNJWat{{\it The weak and strong Lefschetz properties for 
Artinian K-algebras}, J. Algebra {\bf 262} (2003), 99 -- 126.}

\def\MZan{{\it The strength of the weak Lefschetz property},
Illinois J. Math. {\bf 52} (2008), 1417 -- 1433.}

\def\AChoP{{\it Generic initial ideals of artinian ideals having 
Lefschetz properties or the strong Stanley property}, J. Algebra 
{\bf 318} (2007), 589 -- 606.}

\def\BEise{{\it Algebra structures for finite free resolutions
and some structure theorems for ideals of codimension 3}, 
Amer. J. Math. {\bf 99} (1977), 447 -- 485.}

\def\JWat{{\it A note on complete intersections of height three},
Proc. Amer. Math. Soc. {\bf 126} (1998), 3161 -- 3168.}

\makeatletter
\def\@date{September 9, 2011}
\makeatother
\begin{document}

\SkipAmount
\allowdisplaybreaks
\NoFootnoteMark

\RHeads{WLP and basic sequences}
\title{The weak Lefschetz property for Artinian \\ 
graded rings and basic sequences}
\author{\Watashidesu}
\maketitle
\footnotetext{
2010 Mathematics Subject Classification. 
Primary 13D02, 13H10; Secondary 13P10.
}

%
%

\begin{abstract}
The basic sequence of a homogeneous ideal $I\sset R=k[\seq{x}{1}{r}]$ 
defining an Artinian graded ring $A=R/I$ not 
having the weak Lefschetz property has the property that 
the first term of the last part is less than the last term of the penultimate part.
For a general linear form $\ell$ in $\seq{x}{1}{r}$, 
this fact affects in a certain way the behavior of the $r-1$ square 
matrices in $k[\ell]$ which represent the multiplications of the elements of $A$ by
$\seq{x}{1}{r-1}$ through a minimal free presentation of $A$ over $k[\ell]$.
Taking advantage of it, we consider some modules over an algebra generated over 
$k[\ell]$ by the square matrices mentioned above. In this manner, for the case $r=3$, 
we prove that an Artinian \Gor\ graded ring $A=k[x_1,x_2,x_3]/I$ has 
the weak Lefschetz property if $\ch{k}=0$ and the number of the minimal generators
of $0:_A\ell$ over $k[x_1,x_2,x_3]$ is two.
\end{abstract}

\section*{Introduction}
\lb{WLP}

\par
From the viewpoint of standard free resolutions (see \cite[Section 3]{A8}),
homogeneous ideals defining curves in $\Pthree$ and homogeneous ideals
in a polynomial ring $k[x_1,x_2,x_3]$ defining Artinian 
graded rings can be treated in the same manner to a certain extent.
We have recently found a way to apply the method developed in
\cite{A11} to the study of the weak Lefschetz property for Artinian
\Gor\ graded rings. 

\par
In \cite{A11} the structure of the standard free resolutions itself was 
the main theme and the results were expressed in terms of basic sequences
(see \cite[Section 2]{A5}, \cite[Section 1]{A7} or 
\cite[Section 1]{A8} for the definition of basic sequence). 
In this paper we examine the behavior of a ring not having the 
weak Lefschetz property by analyzing a part of its structure which is 
affected by the basic sequence of its defining ideal. 
In contrast to \cite{A11}, 
our arguments are carried out in a way that does not depend too much
on computations of the relation matrices appearing in the standard free 
resolutions, though most of our ideas stem out of observations of them. 
In that manner we can prove that 
an Artinian \Gor\ graded ring $A:=k[x_1,x_2,x_3]/I$ has the weak Lefschetz 
property if $\ch{k}=0$ and the number of the minimal generators of 
$0:_A\ell$ over $k[x_1,x_2,x_3]$ is two for a general linear form 
$\ell$ (see \Thm{\ref{proc61}}). Our results include the case
of height three Artinian complete intersections treated
in \cite[\Thm{2.3}]{HMNJW} and the case where
$\min\{\ l\ |\ \gradn{I}{l}\neq0\ \}=2$ treated in \cite[\Cor{3.2}]{MZ}.

\par
Let us explain the main points shortly.
Let $R:=k[\seq{x}{1}{r}]$ be a polynomial ring over an infinite field
$k$ and $I$ a homogeneous ideal in $R$ such that $A:=R/I$ is an
Artinian graded ring. We assume that $\seq{x}{1}{r}$ are chosen 
sufficiently generally.
In our proof, we consider a minimal free resolution
\begin{equation}
0\lra F\xra{\ \Phg\ }L\xra{\ \pg\ }A\lra 0
\lb{eq52}
\end{equation}
of $A$ over $\Rstar:=k[x_r]$ with homogeneous homomorphisms $\Phg$ 
and $\pg$ of degree zero, where
$F:=\dss{i=1}{m_r}\Rstar(-n^r_i)$,
$L:=\dss{i=1}{m_r}\Rstar(\kg_i)$, 
$(n^r_1\ddd n^r_{m_r})$ is the last part of the basic sequence
$(a;n^2_1\ddd n^2_a;\cdots;n^r_1\ddd n^r_{m_r})$ of $I$
(see \cite[Section 2]{A5}), and $\seq{\kg}{1}{m_r}$ is a nondecreasing
sequence of integers.
For each $\ee{1}{j}{r-1}$, let $\Xsub{j}$ be a matrix  with 
components in $\Rstar$ giving a homogeneous linear map
$\Xsub{j}:L(-1)\lra L$
of degree zero which represents the linear map 
$\times x_j:A(-1)\lra A$ over $\Rstar$. Then $L$ and
$\Ima{\Phg}$ are modules over the $\Rstar$ algebra
$\Rstar[\Xsub{1}\ddd\Xsub{r-1}]$ generated by 
$\Xsub{1}\ddd\Xsub{r-1}$. We see that $A$ has the
weak Lefschetz property if and only if 
$n^r_{m_r}\geq n^{r-1}_{m_{r-1}}$ by \Lem{\ref{proc5}}
and that something special happens to the 
$\Rstar[\Xsub{1}\ddd\Xsub{r-1}]$ module $L$
by \Thm{\ref{proc27}} if $n^r_{m_r} < n^{r-1}_{m_{r-1}}$.
In the case where $A$ is \Gor\ of dimension zero, 
$\ch{k}=0$, $r=3$, and the number of the minimal generators of 
$0:_A\ell$ over $R$ is two, we are led to a contradiction
by the above two facts if $A$ does not have the weak Lefschetz property,
estimating the length of a certain submodule of $\Hom{k}{A/x_3A}{k}
=\Ext{2}{R/(x_3)}{A/x_3A}{R/(x_3)}(-2)$ in two ways, one by
using \eqref{eq52} tensored with $R/(x_3)$ and the other by
using the minimal free resolution of $A/x_3A$ over $R/(x_3)$.
See \Thm{\ref{proc61}} for the detail.

\par
This paper is organized as follows.

\par
In the first two sections, we consider an $R$ module $E$ which is 
finitely generated over $\Rstar:=k[\seq{x}{\rg+1}{r}]$ $(\en{1}{\rg}{r})$ 
with focus on the degrees of its minimal generators over $\Rstar$.
Let 
\begin{equation*}
\cdots\lra F\xra{\ \Phg\ }L\xra{\ \pg\ }E\lra 0
\end{equation*}
be a minimal free resolution of $E$ over $\Rstar$, where
$F:=\dss{i=1}{Q'}\Rstar(-\kg'_i)$, $L:=\dss{i=1}{Q}\Rstar(\kg_i)$,
and $\seq{\kg}{1}{Q}$ and $\kg'_1\ddd\kg'_{Q'}$ are nondecreasing 
sequences of integers. Let further
$\Xsub{j}$ be a matrix  with 
components in $\Rstar$ giving a homogeneous linear map
$\Xsub{j}:L(-1)\lra L$
of degree zero representing the linear map 
$\times x_j:E(-1)\lra E$ over $\Rstar$ for each $\ee{1}{j}{\rg}$.
Under this setting, we can apply the method developed in \cite[Section 3]{A11},
which makes use of the fact that the above sequences of 
integers do not vary when a
small change of variables $\seq{x}{1}{r}$ is performed.
One of our major results here is \Thm{\ref{proc27}} which says that,
for every $p$ satisfying $\kg_p<\kg_{p+1}$, we have
$\Wsup{\ng}\gradn{\Ima{\Phg}}{\leq -\kg_p}\sset N_p$
for all $\ng\geq0$, where 
$W:=\Xsub{1}+\smm{j=2}{\rg}s_j\Xsub{j}$
($\seq{s}{2}{\rg}\in k$), 
$\gradn{\Ima{\Phg}}{\leq -\kg_p}
:=\sm{d\leq -\kg_p}\gradn{\Ima{\Phg}}{d}$, and 
$N_p$ denotes the submodule of $L$ consisting of all of its elements
whose first $p$ components are zero. 
This theorem and the results derived from it play crucial roles 
in the whole paper.

\par
The basic sequence of a module is the sequence of the degrees, 
lined up in a certain rule,
of generators forming a \Wei\ basis of the given module with 
respect to sufficiently general $\seq{x}{1}{r}$ 
(see \cite[Section 2]{A5}). 
It is 
therefore natural that \Wei\ bases appear in our argument frequently.
But, in some cases the conditions imposed on a \Wei\ basis are too
restrictive. For this reason, we will use in this paper a basis satisfying 
relaxed conditions, which we shall call a pseudo \Wei\ basis.
Section three is devoted to a description of its fundamental properties.

\par
In section four some results on the exact 
sequences \eqref{eq52} for Artinian \Gor\ graded rings are described in detail.
At the same time we look into the structure of 
the $\Rstar[\Xsub{1}\ddd\Xsub{r-1}]$ modules $L$ and
$\Ima{\Phg}$ a little.

\par
The main theorem of this paper is proved in section five 
(see \Thm{\ref{proc61}}). There, we work over an Artinian 
\Gor\ graded ring $A=R/I$
such that $r=3$ and the number of the minimal generators of 
$0:_Ax_3$ over $R$ is two, assuming $\ch{k}=0$. 
We first investigate the properties of
the minimal free resolution of $I+(x_3)/(x_3)$ over $R/(x_3)$ and 
that of $I$ over $R$. Then, combining the results obtained in section four, we reach our goal.

\par
In the final section, the basic sequence of a complete intersection of
three homogeneous polynomials is described.

\section{\mathversion{bold}A property coming from invariance under small 
homogeneous transformations of variables}
\lb{kiwadoi}

\par
The principal idea of the arguments of this section
has already been presented in \cite[Section 3]{A11} for
studying a special case. We will enhance it to a more
general formulation with some modifications.

\par
Throughout this paper, let $\seq{y}{1}{r}$ denote 
indeterminates over an infinite field $k$, 
$R$ the polynomial ring 
$k[\seq{y}{1}{r}]$, $\mfrak$ the maximal
ideal $(\seq{y}{1}{r})$,
$\gg_{ij}\ (\ee{1}{i}{r},\ \ee{1}{j}{r})$ elements of
$k$ such that the matrix $\Gg:=(\gg_{ij})$ is invertible, and $\seq{x}{1}{r}$ 
elements of $R$ satisfying
$y_i=\smm{j=1}{r}\gg_{ji}x_j$ $(\ee{1}{i}{r})$. 
Given a matrix $\Xg$ in $R$, the set consisting of all the linear combinations 
of the columns of $\Xg$ over a ring $\mathit{Ring}$ will be denoted by
$\Imasup{\mathit{Ring}}{\Xg}$.
For modules $E,\ E'$ and $E''$, 
the symbol $\op$ will be used in the following two senses\,: (i) $E' \op
E'' = \{\ (e',e'') \ |\ e' \in E', \ e'' \in E'' \ \}$, \ (ii)  $E' \sset E, \ E''
\sset E, \ E' \cap E'' = 0, \ E' \op E'' = \{\ e' + e'' \ |\ e' \in  E', \ e'' 
\in E'' \ \} \sset E$.  Usually the context will make it clear 
which it means.
But in some cases, we will use another symbol $\opdot$ instead of
$\op$ to express the direct sum in the first sense.
The $(Q,Q)$ unit matrix will be denoted by $1_Q$ and   
$\MATr(\Rstar)$ will denote the set of matrices with entries in 
$\Rstar$.

\par
Let $\seq{\kg}{1}{Q}$ be a nondecreasing sequence of integers, 
$\rg$ an integer with $\en{1}{\rg}{r}$, 
and $\Thg_j\ (\ee{1}{j}{\rg})$ matrices giving homogeneous homomorphisms 
\begin{equation*}
\Thg_j:\dss{i=1}{Q}R(\kg_i-1)\lra\dss{i=1}{Q}R(\kg_i)
\end{equation*}
of degree zero 
such that the components of $\Xsub{i}:=x_i1_Q-\Thg_i$ 
lie in $\Rstar:=k[\seq{x}{\rg+1}{r}]$.
Suppose there are positive integers $p,\ q$ with $\en{1}{p}{Q},\ \ee{1}{q}{Q-p}$
satisfying $\kg_p<\kg_{p+1}$ and 
$q=\max\{\ i\ |\ \kg_{p+1}=\cdots=\kg_{p+i},\ p+i\leq Q\ \}$.
Take elements $s_j\in k\ (\ee{2}{j}{\rg})$ and put
$z_1:=x_1+\smm{j=2}{\rg}s_jx_j$, 
$W:=\Xsub{1}+\smm{j=2}{\rg}s_j\Xsub{j}$,
$\Og:=z_11_Q-W=\Thg_1+\smm{j=2}{\rg}s_j\Thg_j$. 

\par
Let
\begin{equation*}
\Og=:\begin{bmatrix}
C'&C''&C'''\\
D'&D''&D'''
\end{bmatrix},\quad 
U':=\begin{bmatrix}
C''&C'''\\
D''&D'''
\end{bmatrix},\quad
U'':=\begin{bmatrix}
C''\\
D''
\end{bmatrix},
\end{equation*}
where the number of the rows of $(C'\ C''\ C''')$
 (resp. $(D'\ D''\ D''')$) is $p$
(resp. $Q-p$), and the numbers of the columns of $D',D''$ and 
$D'''$ are $p,q$ and $Q-p-q$ respectively.
Notice that for each row of $U''$ the degrees of its components are the same
and that $C''$ 
is a matrix whose components are zero or of degree zero in 
$z_1,\seq{x}{\rg+1}{r}$.
Besides, $C'''=0$ since the degrees its components must be negative.
Choose $V_1\in \GL(q,k)$ so that
the columns of $C''V_1$  
different from zero are linearly independent over $k$.

\par
Put
\begin{equation*}
V_2:=
\begin{bmatrix}
1_{p}&0&0\vspace{3pt}\\
0&V_1&0\vspace{3pt}\\
0&0&1_{Q-p-q}
\end{bmatrix}\quad\andt\quad
U:=U'
\begin{bmatrix}
V_1&0\\
0&1_{Q-p-q}
\end{bmatrix}.
\end{equation*}
Observe that $z_1$ appears in $U$ only in the form $cz_1$ with
some $c\in k$. We may therefore write 
\begin{equation}
U=z_1\Usub{1}+\Usub{0},
\lb{eq1}
\end{equation}
where $\Usub{1}$ (resp. $\Usub{0}$) is a matrix with entries 
in $k$ (resp. $\Rstar$).
Let $S$ 
denote the graded module 
$\dss{i=1}{Q}R(\kg_i)$ and 
$\deg(v)$ the degree of an element $v\in S$.
We will regard the columns of $U$ as 
homogeneous elements of $S$.
Note that the degrees of the first $q$ columns of 
$U$ are the same and equal to $-\kg_{p+1}+1$, while the degrees of the remaining
columns are smaller than that.
Let $\seq{b}{1}{m}$ be all the columns of $U$ of degree $-\kg_{p+1}+1$
which do not vanish modulo $(z_1,\seq{x}{\rg+1}{r})$, and denote the 
remaining columns of $U$ by $\seq{a}{1}{n}$, where $n:=Q-p-m$.
Actually, $\{\seq{b}{1}{m}\}$ consists of all the
columns $b$ of $U''V_1$ such that at least one of the first $p$ components 
of $b$ is an element of $k$ different from zero by the choice of $V_1$. 
On the other hand, $\deg(a_i)\leq -\kg_{p+1}+1\leq -\kg_p$ for all 
$\ee{1}{i}{n}$ and the first $p$ components of $a_i$ 
$(\ee{1}{i}{n})$ are zero.
Notice that all the columns of $U$ of degree $-\kg_{p+1}+1$
vanish modulo $(z_1,\seq{x}{\rg+1}{r})$ when $\kg_{p+1}>\kg_p+1$.
In that case, $m=0$ and $n=Q-p$.
Let $b_i=:\tpose{(b_{1i}\ddd b_{Qi})}$ and $b''_i:=\tpose{(b_{1i}\ddd b_{pi})}$.
The components of 
$b''_i$ lie in $k$ and 
the vectors $b''_1\ddd b''_m$ are linearly independent over $k$ by the choice of $V_1$.

\par
We construct a matrix $H$ whose columns are  
$z_1^{\ng_1}x_{\rg+1}^{\ng_{\rg+1}}\cdots x_r^{\ng_r}a_i$ 
$(\ng_1+\smm{j=\rg+1}{r}\ng_j+\deg(a_i)=-\kg_p,\ \ee{1}{i}{n},\ 
\ng_1,\seq{\ng}{\rg+1}{r}\geq0)$ arranged in a suitable order.
Observe that the first $p$ rows of $H$ are zero.
Since $\Og-z_11_Q\in\MATr(\Rstar)$, 
we see $U'-\bclmsc{0}{z_11_{Q-p}}\in\MATr(\Rstar)$,
so that
\begin{equation*}
V_2^{-1}U-
\begin{bmatrix}
O\\
z_11_{Q-p}
\end{bmatrix}
\in \MATr(\Rstar).
\end{equation*}
This implies that
\begin{equation*}
\dss{i=1}{Q}\Rstar[z_1](\kg_i)=
\Imasup{\Rstar[z_1]}{V_2^{-1}U}\op
\lt(\lt(\dss{i=1}{p}\Rstar[z_1](\kg_i)\rt)\opdot
\lt(\dss{i=p+1}{Q}\Rstar(\kg_i)\rt)\rt)
\end{equation*}
by \cite[\Lem{1.1}]{A5}, where 
$\Rstar[z_1]=k[z_1,\seq{x}{\rg+1}{r}]$.
Multiplying both sides of the above equality by 
$V_2$ on the left, we get
\begin{equation*}
\dss{i=1}{Q}\Rstar[z_1](\kg_i)=
\Imasup{\Rstar[z_1]}{U}\op
\lt(\lt(\dss{i=1}{p}\Rstar[z_1](\kg_i)\rt)\opdot
\lt(\dss{i=p+1}{Q}\Rstar(\kg_i)\rt)\rt).
\end{equation*}
Moreover, since an element of 
$\gradn{\Rstar[z_1](\kg_i)}{-\kg_p}$ is zero or lies in 
$k$ for all $\ee{1}{i}{p}$, we see
\begin{equation}
\begin{split}
\gradn{\dss{i=1}{Q}\Rstar[z_1](\kg_i)}{-\kg_p}
&=\gradn{\Imasup{\Rstar[z_1]}{U}}{-\kg_p}\op
\gradn{\dss{i=1}{Q}\Rstar(\kg_i)}{-\kg_p}\\
&=\la \seq{b}{1}{m},\ H\ra \op
\gradn{\dss{i=1}{Q}\Rstar(\kg_i)}{-\kg_p},
\end{split}
\lb{eq2}
\end{equation}
where $\la \seq{b}{1}{m},\ H\ra$ denotes the vector space over $k$
spanned by $\seq{b}{1}{m}$ and the columns of $H$.
Note that $\la \seq{b}{1}{m},\ H\ra=\la H\ra$ when $\kg_{p+1}>\kg_p+1$.

\par
We will denote a matrix $Z$ with components
in $R$ by $Z(\seq{x}{1}{r})$ when we want to pay attention to
the variables $\seq{x}{1}{r}$.
Let $\xg_j\ (\ee{\rg+1}{j}{r})$ be parameters over $R$.
For a matrix $Z=Z(\seq{x}{1}{r})$ with components in $R$,
let  
\begin{align*}
\Ztild=\Ztild(\seq{\xg}{\rg+1}{r},\seq{x}{1}{r})
:&=Z(\seq{x}{1}{\rg},x_{\rg+1}+\xg_{\rg+1}z_1\ddd x_r+\xg_rz_1)\\
\intertext{and} 
\Zbar=\Zbar(\seq{\xg}{\rg+1}{r},\seq{x}{1}{\rg})
:&=Z(\seq{x}{1}{\rg},\xg_{\rg+1}z_1\ddd\xg_rz_1)\\
&=\Ztild(\seq{\xg}{\rg+1}{r},\seq{x}{1}{\rg},0\ddd 0).
\end{align*}
Observe that 
\begin{equation}
\Zbar(x_{\rg+1}/z_1\ddd x_r/z_1,\seq{x}{1}{\rg})=Z(\seq{x}{1}{r}).
\lb{eq3}
\end{equation}
Now we consider $\seq{b}{1}{m}$ and $H$. Notice that their components are
polynomials in $z_1,\seq{x}{\rg+1}{r}$. 
With the notation above, the components of 
$\seq{\btild}{1}{m}$ and $\Htild$, therefore, lie in 
$k[z_1,\seq{x}{\rg+1}{r},\seq{\xg}{\rg+1}{r}]
=\Rstar[z_1]\ot_k k[\seq{\xg}{\rg+1}{r}]$. 
Let $T$ be the local ring $k[\seq{\xg}{\rg+1}{r}]_{(\seq{\xg}{\rg+1}{r})}$ and
let $\la \seq{\btild}{1}{m},\ \Htild\ra_T$ denote 
the submodule of 
\begin{equation*}
\gradn{\dss{i=1}{Q}\Rstar[z_1](\kg_i)
\ot_k T}{-\kg_p}
\end{equation*}
spanned over $T$ by $\seq{\btild}{1}{m}$ and the columns of $\Htild$,
where $\gradn{S\ot_k T}{j}=\gradn{S}{j}\ot_k T$ for $j\in\Zbf$.
Note that $\la \seq{\btild}{1}{m},\ \Htild\ra_T=\la\Htild\ra_T$ in the case 
$\kg_{p+1}>\kg_p+1$.
Since $\btild_i(0\ddd 0,\seq{x}{1}{r})=b_i\ (\ee{1}{i}{m})$ and 
$\Htild(0\ddd 0,\seq{x}{1}{r})=H$, we find by
\eqref{eq2} that 
\begin{equation}
\gradn{\dss{i=1}{Q}\Rstar[z_1](\kg_i)
\ot_k T}{-\kg_p}
=\la \seq{\btild}{1}{m},\ \Htild\ra_T \op
\gradn{\dss{i=1}{Q}\Rstar(\kg_i)\ot_k T}{-\kg_p}.
\lb{eq4}
\end{equation}
Put
\begin{equation*}
N_p:=0^p\opdot 
\lt(\dss{i=p+1}{Q}\Rstar(\kg_i)\rt)
\sset\lt(\dss{i=1}{Q}\Rstar(\kg_i)\rt).
\end{equation*}
This is a module over $\Rstar$ consisting of all the elements of 
$\dss{i=1}{Q}\Rstar(\kg_i)\sset\dss{i=1}{Q}R(\kg_i)=S$
such that the first $p$ components are zero. 
Let $v$ be an element of 
$\grad{N_p}{-\kg_p}$. Note that no variables $\seq{x}{1}{\rg}$
appear in $v$ although we write $v=v(\seq{x}{1}{r})$.
We have
\begin{align}
&\vtild(\seq{\xg}{\rg+1}{r},\seq{x}{1}{r})
\notag\\
&=
\smm{i=1}{m}\btild_i(\seq{\xg}{\rg+1}{r},\seq{x}{1}{r})\ghat_i+
\Htild(\seq{\xg}{\rg+1}{r},\seq{x}{1}{r})\cdot\tpose{\row{\fhat}{1}{l}}
\notag\\
 &\hspace{8em} +w
\lb{eq5}
\end{align}
with $\ghat_i,\ \fhat_i\in T$ and 
$w\in \gradn{\dss{i=1}{Q}\Rstar(\kg_i)\ot_k T}{-\kg_p}$ by \eqref{eq4},
where $l$ denotes the number of the columns of $H$.

\begin{lem}\lb{proc2}
Let $v$ and $w$ be as above. Suppose that 
$w\equiv 0\ (\mod\ (\seq{x}{\rg+1}{r}))$. Then 
\begin{equation*}
z_1^\ng v\in \Imasup{\Rstar[z_1]}{\seq{a}{1}{n}}
\op N_p
\end{equation*}
for all $\ng\geq0$.
\end{lem}

\begin{pf}
Since $v\in\gradn{N_p}{-\kg_p}$ by hypotheses, 
the first $p$ components of $\vtild$
are zero. Besides, since an element of 
$\gradn{\Rstar[z_1](\kg_i)\ot_k T}{-\kg_p}$ 
is zero or lies in
$T$ for $\ee{1}{i}{p}$, it follows from the assumption 
$w\equiv 0\ (\mod\ (\seq{x}{\rg+1}{r}))$
that the first $p$ components of $w$ are also zero.
On the other hand, the vectors $\btild''_i$
$(\ee{1}{i}{m})$ are linearly independent over $T$ and the first $p$ 
components of the columns of
$\Htild$ are zero, since these properties
are inherited from $b''_1\ddd b''_m$ and $H$.
Hence, $\ghat_i=0$ for all $\ee{1}{i}{m}$. In other words
\begin{equation*}
\vtild(\seq{\xg}{\rg+1}{r},\seq{x}{1}{r})
=\Htild(\seq{\xg}{\rg+1}{r},\seq{x}{1}{r})
\cdot\tpose{\row{\fhat}{1}{l}}+w,
\end{equation*}
so that
\begin{equation*}
\vbar(\seq{\xg}{\rg+1}{r},\seq{x}{1}{\rg})
=\Hbar(\seq{\xg}{\rg+1}{r},\seq{x}{1}{\rg})
\cdot\tpose{\row{\fhat}{1}{l}}.
\end{equation*} 
Let $\ng\geq0$ be an integer. Since the denominators of $\fhat_i\ (\ee{1}{i}{l})$ lies 
in $k^\ast+(\seq{\xg}{\rg+1}{r})$, there is a polynomial 
$\psg_0\in(\seq{\xg}{\rg+1}{r})^{\ng+1}$
such that $\psg_i:=(1+\psg_0)\fhat_i\in k[\seq{\xg}{\rg+1}{r}]$ 
for all $\ee{1}{i}{l}$. Hence
\begin{equation*}
(1+\psg_0)z_1^\ng\vbar(\seq{\xg}{\rg+1}{r},\seq{x}{1}{\rg})
=z_1^\ng\Hbar(\seq{\xg}{\rg+1}{r},\seq{x}{1}{\rg})
\cdot\tpose{\row{\psg}{1}{l}}.
\end{equation*} 
Now substitute $x_j/z_1$ for $\xg_j$ for all $\ee{\rg+1}{j}{r}$
in this equality.
We find by \eqref{eq3} that
\begin{multline}
(1+\psg_0(x_{\rg+1}/z_1\ddd x_r/z_1))z_1^\ng 
v(\seq{x}{\rg+1}{r})\\=
z_1^\ng H(\seq{x}{1}{r})
\cdot\tpose{(\psg_1(x_{\rg+1}/z_1\ddd x_r/z_1),
\ldots,\psg_l(x_{\rg+1}/z_1\ddd x_r/z_1))}.
\lb{eq19}
\end{multline}
Write
\begin{equation*}
\psg_j(x_{\rg+1}/z_1\ddd x_r/z_1)=
\sm{\mg\geq0}\psg_{j\mg}(\seq{x}{\rg+1}{r})/z_1^\mg
\end{equation*}
for $\ee{0}{j}{l}$,
where $\psg_{j\mg}(\seq{x}{\rg+1}{r})$ is a homogeneous polynomial in 
$\seq{x}{\rg+1}{r}$ of degree $\mg$ for each $j,\mg$. 
Notice that $\psg_{0\mg}=0$ for
$\mg\leq\ng$. Moreover, we can write 
$a_i=z_1a_{i1}+a_{i0}$ by \eqref{eq1}, where $a_{i1}$ (resp. $a_{i0}$)
is a column of $\Usub{1}$ (resp. $\Usub{0}$).
Compare terms with no factor $z_1$ in the denominators in the 
above equality  \eqref{eq19}.
Then, since $v\in\dss{i=1}{Q}\Rstar(\kg_i)$ and each column  of 
$H$ is of the form
$z_1^{\ng_1}x_{\rg+1}^{\ng_{\rg+1}}\ddd x_r^{\ng_r}a_i$, we find that
$z_1^\ng v$ is the sum of a finite number of vectors of the forms
\begin{align*}
&z_1^{\ng_1+\ng-\mg}x_{\rg+1}^{\ng_{\rg+1}}\ddd x_r^{\ng_r}
\psg_{j\mg}(\seq{x}{\rg+1}{r}) a_i \\
&\hskip3em(\deg(a_i)+\ng_1+\smm{t=\rg+1}{r}\ng_t=-\kg_p,\ 
\ng_t\geq0\ (t=1,\rg+1\ddd r),\\
&\hskip6em\ng_1+\ng\geq\mg,\
\ee{1}{i}{n},\ \ee{1}{j}{l})\quad \andr\\
&x_{\rg+1}^{\ng_{\rg+1}}\ddd x_r^{\ng_r}
\psg_{j\ng_1+\ng+1}(\seq{x}{\rg+1}{r}) a_{i1} \\
&\hskip3em(\deg(a_i)+\ng_1+\smm{t=\rg+1}{r}\ng_t=-\kg_p,\ 
\ng_t\geq0\ (t=1,\rg+1\ddd r),\\ 
&\hskip6em \ee{1}{i}{n},\ \ee{1}{j}{l}).
\end{align*}
Since the first $p$ components of $a_{i1}$
are zero for all $\ee{1}{i}{n}$, our assertion holds. 
\end{pf}

\begin{lem}\lb{proc3}
Assumption being the same as in \Lem{\ref{proc2}}, we have
\begin{equation*}
\Wsup{\ng} v\in N_p
\end{equation*}
for all $\ng\geq0$.
\end{lem}

\begin{pf}
For $\ng=0$, our assertion is trivial. Suppose $\ng>0$.
Using the equality $\Og=z_1 1_Q-W$, we see 
\begin{equation*}
z_1^\ng v
=\Og\lt(\smm{i=1}{\ng}z_1^{\ng-i}\Wsup{i-1}\rt)v
+\Wsup{\ng} v.
\end{equation*}
On the other hand, $z_1^\ng v=\og_1+\og_2$ with 
$\og_1\in\Imasup{\Rstar[z_1]}{\seq{a}{1}{n}}
\sset \Imasup{\Rstar[z_1]}{\Og}$, 
$\og_2\in N_p$ by \Lem{\ref{proc2}}.
The vector $\og_2-\Wsup{\ng} v$ is therefore contained in
$\Imasup{\Rstar[z_1]}{\Og}$. Since the components of 
$\og_2-\Wsup{\ng} v$ are all contained in $\Rstar$ and 
$\Og$ is a
square matrix such that $\Og=z_1 1_Q-W$ with
$W\in \MATr(\Rstar)$,
it follows that $\og_2-\Wsup{\ng} v=0$. Hence $\Wsup{\ng} v\in N_p$.
\end{pf}

\begin{rem}\lb{proc34}
When $\kg_{p+1}>\kg_p+1$, the $j$-th column of $W$ lies in
$N_p$ for all $j$ with $\ee{p+1}{j}{Q}$. Hence $\Wsup{\ng}N_p\sset N_p$
for all $\ng\geq0$.
\end{rem}

\begin{lem}\lb{proc28}
With the notation above, 
\begin{equation*}
S=\dss{j=1}{\rg}\Imasup{k[\seq{x}{j}{r}]}{\Thg_j}\op
\dss{i=1}{Q}\Rstar(\kg_i).
\end{equation*}
\end{lem}

\begin{pf}
By the same argument as that of the proof of \cite[\Lem{1.1}]{A5},
\begin{equation*}
\dss{i=1}{Q}k[\seq{x}{j}{r}](\kg_i)
=\Imasup{k[\seq{x}{j}{r}]}{\Thg_j}\op
\dss{i=1}{Q}k[\seq{x}{j+1}{r}](\kg_i)
\end{equation*}
for all $\ee{1}{j}{\rg}$. Our assertion follows by repeated use of
this formula.
\end{pf}

\begin{lem}\lb{proc4}
Let $\kg_i$ $(\ee{1}{i}{Q})$, $\Thg_j$ $(\ee{1}{j}{\rg})$, 
$p$, $q$, $N_p$, $\xg_j$ $(\ee{\rg+1}{j}{r})$, $S$ and
$T$ be as above. 
Let $\seq{v}{1}{\dg}$ be homogeneous elements of 
$\dss{i=1}{Q}\Rstar(\kg_i)$ and 
$v$ an element of 
$\gradn{N_p}{-\kg_p}\cap\Imasup{\Rstar}{\seq{v}{1}{\dg}}$.
Suppose that the number of the homogeneous minimal 
generators of 
\begin{equation*}
S/\Imasup{R}{\Thg_1\ddd\Thg_\rg,\seq{v}{1}{\dg}}
\end{equation*}
over $\Rstar$ of degree $-\kg_p$ is 
$q':=\max\{\ i\ |\ \kg_{p-i+1}=\cdots=\kg_{p-1}=\kg_p\ \}$ 
and remains unchanged for any small
homogeneous transformation of variables 
$\seq{x}{1}{r}$.
Then, for every sequence $s_2\ddd s_\rg\in k$, the element $w$ of 
\begin{equation*}
\gradn{\dss{i=1}{Q}\Rstar(\kg_i)\ot_k T}{-\kg_p}
\end{equation*}
defined by the equality
\eqref{eq5} must be congruent to zero modulo $(\seq{x}{\rg+1}{r})$. 
\end{lem}

\begin{pf}
Let $s_2\ddd s_\rg\in k$, and let $z_1$, $\Og$, $U$, 
$b_i\ (\ee{1}{i}{m})$, and $H$ be as above.
We want to consider the $R\ot_k k[\seq{\xg}{\rg+1}{r}]$-module
$\Coksup{R\ot_k k[\seq{\xg}{\rg+1}{r}]}
{\seq{\Thgtild}{1}{\rg},\seq{\vtild}{1}{\dg}}$
as a module over $\Rstar\ot_k k[\seq{\xg}{\rg+1}{r}]$,
where $\Thgtild_i=\Thgtild_i(\seq{\xg}{\rg+1}{r},\seq{x}{1}{r})\ 
(\ee{1}{i}{\rg})$ and 
$\vtild_j=\vtild_j(\seq{\xg}{\rg+1}{r},\seq{x}{1}{r})$ $(\ee{1}{j}{\dg})$.
We see 
\begin{equation*}
S=\Imasup{R}{\seq{\Thg}{1}{\rg}}
+\dss{i=1}{Q}\Rstar(\kg_i)
\end{equation*}
by \Lem{\ref{proc28}}, so that
\begin{equation*}
S\ot_k T=\Imasup{R\ot_k T}{\seq{\Thgtild}{1}{\rg}}+
\dss{i=1}{Q}\Rstar(\kg_i)\ot_k T.
\end{equation*}
Hence 
\begin{align}
&\lt(S\ot_k k[\seq{\xg}{\rg+1}{r}]_\psg\rt)/
\Imasup{R\ot_k k[\seq{\xg}{\rg+1}{r}]_\psg}
{\seq{\Thgtild}{1}{\rg},\seq{\vtild}{1}{\dg}}
\notag\\
&\hspace{4em}\cong 
\frac{\dss{i=1}{Q}\Rstar(\kg_i)\ot_k 
k[\seq{\xg}{\rg+1}{r}]_\psg}
{(\Rstar\ot_k k[\seq{\xg}{\rg+1}{r}]_\psg)E'}
\lb{eq6}
\end{align}
over $\Rstar\ot_k k[\seq{\xg}{\rg+1}{r}]_\psg$ for some 
$\psg\in k[\seq{\xg}{\rg+1}{r}]\bslash(\seq{\xg}{\rg+1}{r})$, where
\begin{equation*}
\begin{split}
E':=&\Imasup{R\ot_k k[\seq{\xg}{\rg+1}{r}]}
{\seq{\Thgtild}{\rg+1}{r},\seq{\vtild}{1}{\dg}}\\
&\cap \lt(\dss{i=1}{Q}\Rstar(\kg_i)\ot_k k[\seq{\xg}{\rg+1}{r}]\rt).
\end{split}
\end{equation*}
Let 
$e_i:=\tpose{(0\ddd0,\overset{\indexonsmile{i}}1,0\ddd0)}$
$(\ee{1}{i}{Q})$ denote the canonical bases of 
$\dss{i=1}{Q}\Rstar(\kg_i)\ot_k k[\seq{\xg}{\rg+1}{r}]$. 
Then $\deg(e_i)=-\kg_p$
if and only if $\ee{p-q'+1}{i}{p}$. Since the parameters 
$\seq{\xg}{\rg+1}{r}$ correspond to a small homogeneous transformation of
the variables $\seq{x}{1}{r}$, our hypothesis implies by 
\eqref{eq6} that the $q'$ vectors
$\seq{e}{p-q'+1}{p}$ must be linearly independent over 
$k(\pfrak)$ when considered in 
\begin{equation*}
\lt(\frac{\dss{i=1}{Q}\Rstar(\kg_i)\ot_k 
k[\seq{\xg}{\rg+1}{r}]_{\psg}}
{(\Rstar\ot_k k[\seq{\xg}{\rg+1}{r}]_{\psg})E'}\rt)
\ot k(\pfrak)
\end{equation*}
for all points $\pfrak\in\Spec{k[\seq{\xg}{\rg+1}{r}]}$ in a neighborhood of
the origin $\xg_{\rg+1}=\cdots=\xg_r=0$.
We find therefore that the first $p$ components
of any element of $E'$ of degree $-\kg_p$ must be zero.
Let $w$ be the element of
$\dss{i=1}{Q}\Rstar(\kg_i)\ot_k T$ of 
degree $-\kg_p$ such that the equality \eqref{eq5} holds.
The vectors $\btild_i\ (\ee{1}{i}{m})$ and the columns of $\Htild$  
are contained in
$\Imasup{R\ot_k T}{\seq{\Thgtild}{1}{\rg}}$, and 
$\vtild\in\Imasup{R\ot_k T}{\seq{\vtild}{1}{\dg}}$, so that 
$w\in (\Rstar\ot_k T)E'$.
Hence $w\equiv 0 \ (\mod\ (\seq{x}{\rg+1}{r}))$.
\end{pf}

\section{\mathversion{bold}Matrices that represent
operations of $\seq{x}{1}{\rg}$}
\lb{mat}

\par
Let $E$ be a \fgg\ module over $R$. Suppose that there is an integer 
$\rg$ satisfying $\en{1}{\rg}{r}$ such that $E$ has a minimal free 
resolution of the form
\begin{equation}
\lra\dss{i=1}{Q'}\Rstar(-\kg'_i)\xra{\ \Phg\ }
\dss{i=1}{Q}\Rstar(\kg_i)\xra{\ \pg\ } E\lra 0
\lb{eq7}
\end{equation}
over $\Rstar=k[\seq{x}{\rg+1}{r}]$ with homogeneous homomorphisms
$\Phg$ and $\pg$.
We assume here that the sequences $\kg_1\ddd\kg_Q$ and
$\kg'_1\ddd\kg'_{Q'}$ are nondecreasing. Since $E$ is a module over $R$, 
for each $\ee{1}{j}{\rg}$, there is a matrix $\Xsub{j}$ with 
components in $\Rstar$ giving a homogeneous linear map
\begin{equation*}
\Xsub{j}:\dss{i=1}{Q}\Rstar(\kg_i-1)\lra\dss{i=1}{Q}\Rstar(\kg_i)
\end{equation*}
of degree zero which represents the linear map 
\begin{equation*}
\times x_j:E(-1)\lra E
\end{equation*}
over $\Rstar$. 

\par
Put $\Thg_j:=x_j1_Q-\Xsub{j}$ $(\ee{1}{j}{\rg})$.
Then each matrix $\Thg_j$ gives a
homogenous homomorphism
\begin{equation*}
\Thg_j:\dss{i=1}{Q}R(\kg_i-1)\lra\dss{i=1}{Q}R(\kg_i)
\end{equation*}
of degree zero. 
In what follows, for the sake of simplicity,
let $S$ and $L$ denote the graded modules
$\dss{i=1}{Q}R(\kg_i)$ and 
$\dss{i=1}{Q}\Rstar(\kg_i)$ respectively,
and $e_i$ the canonical base 
$\tpose{(0\ddd0,\overset{\indexonsmile{i}}1,0\ddd0)}$
of $L\sset S$ for $\ee{1}{i}{Q}$.
We define a homogeneous linear map $\pg^R:S\lra E$ over $R$ of degree
zero by setting $\pg^R(e_i)=\pg(e_i)$ for all $\ee{1}{i}{Q}$.
Notice that 
$\pg^R(x_jv)=x_j\pg^R(v)=x_j\pg(v)=
\pg(\Xsub{j}v)=\pg^R(\Xsub{j}v)$
for all $v\in L$ and $\ee{1}{j}{\rg}$.

\par
As in the previous section, suppose there are positive 
integers $p,\ q$ with $\en{1}{p}{Q},\ \ee{1}{q}{Q-p}$
satisfying $\kg_p<\kg_{p+1}$ and 
$q=\max\{\ i\ |\ \kg_{p+1}=\cdots=\kg_{p+i},\ p+i\leq Q\ \}$.
Let $s_j$ $(\ee{2}{j}{\rg})$ be arbitrary elements of $k$, and 
define $z_1$, $W$, $\Og$ and $N_p$ in the same way as before. 
Further we put
\begin{equation*}
\gradn{E}{\leq \kgtild}
:=\sm{d\leq \kgtild}\gradn{E}{d}
\sset E
\end{equation*}
for a graded module $E$ and $\kgtild\in\Zbf$.

\begin{thm}\lb{proc27}
With the notation above, suppose that
$\Gg$ is chosen sufficiently generally.
Then,  
$\Wsup{\ng}\gradn{\Ima{\Phg}}{\leq -\kg_p}\sset N_p$
for all $\ng\geq0$.
\end{thm}

\begin{pf}
Since each column of $\Thg_j$ lies in $\Ker{\pg^R}$ for every 
$\ee{1}{j}{\rg}$, one sees that 
\begin{equation*}
\smm{j=1}{\rg}\Imasup{k[\seq{x}{j}{r}]}{\Thg_j}
\sset\smm{j=1}{\rg}\Imasup{R}{\Thg_j}\sset\Ker{\pg^R}.
\end{equation*}
By \Lem{\ref{proc28}}, therefore, 
\begin{equation*}
\begin{split}
\Ker{\pg^R}&=\smm{j=1}{\rg}\Imasup{k[\seq{x}{j}{r}]}{\Thg_j}
+\Ker{\pg^R}\cap L
=\smm{j=1}{\rg}\Imasup{k[\seq{x}{j}{r}]}{\Thg_j}
+\Ker{\pg}\\
&=\smm{j=1}{\rg}\Imasup{R}{\Thg_j}+\Ima{\Phg}
=\smm{j=1}{\rg}\Imasup{R}{\Thg_j}+\Imasup{R}{\Phg}.
\end{split}
\end{equation*}
Let $\phgsub{i}$ denote the $i$-th column of $\Phg$ for $\ee{1}{i}{Q'}$.
Then 
\begin{equation*}
E=S/\Imasup{R}{\seq{\Thg}{1}{\rg},\, \phgsub{1}\ddd\phgsub{Q'}}
\end{equation*}
by what we have seen. 
Now suppose that $\Gg$ is chosen sufficiently generally.
The nondecreasing
sequences $\kg_1\ddd\kg_Q$ and $\kg'_1\ddd\kg'_{Q'}$ appearing in 
\eqref{eq7} do not change under any small homogeneous transformation 
of the variables $\seq{x}{1}{r}$. 
Let $v$ be an arbitrary element of $\gradn{\Ima{\Phg}}{d}$
for some $d\leq -\kg_p$. Since the resolution \eqref{eq7} is
minimal over $\Rstar$, no element of $k^\ast$ occurs in $\Phg$.
This means that 
$x_r^{-\kg_p-d}v
\in\gradn{N_p}{-\kg_p}\cap\Imasup{\Rstar}{\phgsub{1}\ddd\phgsub{Q'}}$.
Hence 
$x_r^{-\kg_p-d}\Wsup{\ng}v
=\Wsup{\ng}x_r^{-\kg_p-d}v\sset N_p$
for all $\ng\geq0$ by \Lems{\ref{proc3} and \ref{proc4}}.
Thus $\Wsup{\ng}v\sset N_p$ for all $\ng\geq0$.
This proves our assertion.
\end{pf}

\par
For each $\ee{1}{j}{\rg}$, 
one sees that $\Xsub{j}\Ima{\Phg}
\sset\Ima{\Phg}$,
since $\pg(\Xsub{j}\Ima{\Phg})
=x_j\pg(\Ima{\Phg})=0$.
Namely, $\Ima{\Phg}$ is a submodule of $L$ over the $\Rstar$ algebra
\linebreak
$\Rstar[\Xsub{1}\ddd\Xsub{\rg}]$.
Varying $\seq{s}{2}{\rg}$, one can find something about 
$\Rstar[\Xsub{1}\ddd\Xsub{\rg}]v$ for $v\in L$. 
But, it should be done
carefully since $\Rstar[\Xsub{1}\ddd\Xsub{\rg}]$ is not 
commutative in general. 
In the following argument, we consider the columns of 
$[\Xsub{i},\Xsub{j}]=\Xsub{i}\Xsub{j}-\Xsub{j}\Xsub{i}$ as elements of $L$.
Let $\thgsub{ijl}$ denote the $l$-th column of
$[\Xsub{i},\Xsub{j}]$. 

\begin{lem}\lb{proc29}
Let the notation be as above.
\begin{parenuma}
\item\lb{cond6}
For every $\ee{1}{l}{Q}$, we have
$\thgsub{ijl}\in\Ima{\Phg}$.
\item\lb{cond7}
The degree of $\thgsub{ijl}$ is $-\kg_l+2$ for each $\ee{1}{l}{Q}$.
\item\lb{cond8}
$[\Xsub{i},\Xsub{j}]N_p
\sset\Rstar\gradn{\Ima{\Phg}}{\leq -\kg_p}
+\Imasup{\Rstar}{\thgsub{ijp+1}\ddd\thgsub{ijp+q}}$.
\end{parenuma}
\end{lem}

\begin{pf}
Left to the readers.
\end{pf}

\begin{thm}\lb{proc18}
With the notation above, suppose that $\Gg$ is chosen sufficiently generally.
Suppose further that $\ch{k}=0$ and that 
$\thgsub{ijp+1}\ddd\thgsub{ijp+q}$ are contained in the submodule
$\Rstar[\Xsub{1}\ddd\Xsub{\rg}]
\gradn{\Ima{\Phg}}{\leq -\kg_p}$ for all $i,j$ with
$1\leq i<j \leq\rg$. 
Then we have $\Rstar[\Xsub{1}\ddd\Xsub{\rg}]
\gradn{\Ima{\Phg}}{\leq -\kg_p}\sset N_p$.
\end{thm}

\begin{pf}
We can write
\begin{equation*}
W^\ng=\sm{\begin{smallmatrix}
\ng_1+\cdots+\ng_\rg=\ng, \\ \ng_i\geq0\ (\ee{1}{i}{\rg})
\end{smallmatrix}}
s_2^{\ng_2}\cdots s_\rg^{\ng_\rg}\vartheta_{\ng_1\ng_2\ldots\ng_\rg},
\end{equation*}
where $\vartheta_{\ng_1\ng_2\ldots\ng_\rg}$ is the sum of monomials of
the form $\Xsub{i_1}\Xsub{i_2}\cdots\Xsub{i_\ng}$ with
$\#\{\ l\ |\ i_l=j\ \}=\ng_j$ for all $\ee{1}{j}{\rg}$.
Since $\seq{s}{2}{\rg}$ can be chosen arbitrarily, one sees by
\Lem{\ref{proc27}} that
\begin{equation}
\vartheta_{\ng_1\ng_2\ldots\ng_\rg}
\gradn{\Ima{\Phg}}{\leq -\kg_p}\sset N_p.
\lb{eq21}
\end{equation}
To prove our assertion, it is enough to show that 
\begin{equation}
\Xsub{i_1}\Xsub{i_2}\cdots\Xsub{i_\ng}
\gradn{\Ima{\Phg}}{\leq -\kg_p}
\sset N_p
\lb{eq18}
\end{equation}
for all sequences $\seq{i}{1}{\ng}\in\{1\ddd\rg\}$ $(\ng\geq0)$.
The cases $\ng=0,1$ follow directly from \Thm{\ref{proc27}}.
Before proceeding further, notice that
\begin{equation}
[\Xsub{i},\Xsub{j}]N_p 
\sset
\Rstar\gradn{\Ima{\Phg}}{\leq -\kg_p}
+\smm{i=1}{\rg}
\Xsub{i}\gradn{\Ima{\Phg}}{\leq -\kg_p}
\lb{eq20}
\end{equation}
by \Lem{\ref{proc29}} and our hypothesis.
Since $\gradn{\Ima{\Phg}}{\leq -\kg_p}\sset N_p$
and $2\Xsub{i}\Xsub{j}
=\Xsub{i}\Xsub{j}+\Xsub{j}\Xsub{i}
+[\Xsub{i},\Xsub{j}]$, 
the case $\ng=2$ follows from \eqref{eq21}, \eqref{eq20} 
and the case $\ng=1$. 
For $\ng>2$ and for each monomial 
$\Xsub{i_1}\Xsub{i_2}\cdots\Xsub{i_\ng}$ with
$\#\{\ l\ |\ i_l=j\ \}=\ng_j$ for every $\ee{1}{j}{\rg}$,
we find that
\begin{align*}
&\frac{\ng!}{\ng_1!\ng_2!\cdots\ng_\rg!}\,
\Xsub{i_1}\Xsub{i_2}\cdots\Xsub{i_\ng}
\\
&\in \vartheta_{\ng_1\ng_2\ldots\ng_\rg}+
\sm{\begin{smallmatrix}
1\leq i<j \leq\rg,\ \mg\leq\ng-2,\\ 
\ee{1}{j_l}{\rg}\ (\ee{1}{l}{\ng-2})                     
\end{smallmatrix}}
\Rstar
\Xsub{j_1}\Xsub{j_2}\cdots\Xsub{j_\mg}
[\Xsub{i},\Xsub{j}]
\Xsub{j_{\mg+1}}\Xsub{j_{\mg+2}}
\cdots\Xsub{j_{\ng-2}},
\end{align*}
making repeated use of $\Xsub{i'}\Xsub{j'}
=\Xsub{j'}\Xsub{i'}+[\Xsub{i'},\Xsub{j'}]$.
Assuming that \eqref{eq18} is true for all smaller values of $\ng$, we
see
\begin{align*}
&\frac{\ng!}{\ng_1!\ng_2!\cdots\ng_\rg!}\,
\Xsub{i_1}\Xsub{i_2}\cdots\Xsub{i_\ng}v
\\
&\in \vartheta_{\ng_1\ng_2\ldots\ng_\rg}v+
\sm{\begin{smallmatrix}
1\leq i<j \leq\rg,\ \mg\leq\ng-2,\\ 
\ee{1}{j_l}{\rg}\ (\ee{1}{l}{\mg})                     
\end{smallmatrix}}
\Rstar
\Xsub{j_1}\Xsub{j_2}\cdots\Xsub{j_\mg}
[\Xsub{i},\Xsub{j}]N_p.
\end{align*}
for all $v\in\gradn{\Ima{\Phg}}{\leq -\kg_p}$.
This, together with \eqref{eq20}, leads us to our assertion
by induction on $\ng$.
\end{pf}

\begin{cor}\lb{proc13}
Let the notation be as above and suppose that $\Gg$ is chosen 
sufficiently generally.
If $\ch{k}=0$ and
$\gradn{\Ima{\Phg}}{-\kg_{p+1}+2}
\sset\Rstar[\Xsub{1}\ddd\Xsub{\rg}]
\gradn{\Ima{\Phg}}{\leq -\kg_p}$, then 
\linebreak
$\Rstar[\Xsub{1}\ddd\Xsub{\rg}]\gradn{\Ima{\Phg}}{\leq -\kg_p}\sset N_p$.
\end{cor}

\begin{pf}
One has $\thgsub{ijp+1}\ddd\thgsub{ijp+q}
\in\gradn{\Ima{\Phg}}{-\kg_{p+1}+2}$ by
\eqref{cond7} of \Lem{\ref{proc29}}.
\end{pf}

\begin{prop}\lb{proc55}
With the notation above, assume that $\rg=r-1$. 
Suppose that $0:_E x_r$ is generated over $R$ by 
$\gradn{0:_E x_r}{\leq\kgtild-1}$ for an
integer $\kgtild$. Then the following assertions
hold.
\begin{parenuma}
\item\lb{cond25}
$\Ima{\Phg}$ is generated over 
$\Rstar[\Xsub{1}\ddd\Xsub{\rg}]$ by 
$\gradn{\Ima{\Phg}}{\leq\kgtild}$.
\item\lb{cond26}
Suppose further that $\Gg$ is chosen sufficiently generally, 
that $\ch{k}=0$, and that $\kgtild+\kg_1\leq0$. Let $p$ be the largest 
among the integers $i$ satisfying $\kgtild+\kg_i\leq0$, $\ee{1}{i}{Q}$.
Then $\Ima{\Phg}\sset N_p$ and $\len{R}{E}=\infty$.
\item\lb{cond38}
(\Cor{3.11} (Socle Lemma) (ii) of \cite{HUl})\quad 
If $\ch{k}=0$, $\len{R}{E}<\infty$ and $\Gg$ is sufficiently general, 
we have $\kgtild > \max\{\ d\ |\ \gradn{E/x_rE}{d}\neq0\ \}$. 
\end{parenuma}
\end{prop}

\begin{pf}
Put $\Rbarstar:=\Rstar/x_r\Rstar$.
When $\rg=r-1$, the components of $\Phg$ lie in
$x_r\Rstar$. Besides, $\pg(v/x_r)\in 0:_Ex_r$ for all 
$v\in\Ima{\Phg}$. Let $\phi$ denote the map from
$\dss{i=1}{Q'}\Rbarstar(-\kg'_i+1)=\dss{i=1}{Q'}k(-\kg'_i+1)$ to 
$0:_Ex_r$ defined by $\phi(\wbar):=\pg(\Phg(w)/x_r)$ for 
$\wbar\in\dss{i=1}{Q'}k(-\kg'_i+1)$, where 
$w\in\dss{i=1}{Q'}\Rbarstar(-\kg'_i+1)$ and $\wbar=w\ (\mod\ (x_r))$.
It can be verified without difficulty that $\phi$ is a homogeneous isomorphism
of degree zero over $k$. 
\begin{parenuma}
\item
It follows from the hypothesis that 
\begin{equation*}
\Imasup{\Rstar}{\Phg/x_r}\sset
\Rstar[\Xsub{1}\ddd\Xsub{\rg}]
\gradn{\Imasup{\Rstar}{\Phg/x_r}}{\leq\kgtild-1}+
\Ima{\Phg},
\end{equation*}
in other words,
\begin{equation*}
\Ima{\Phg}\sset
\Rstar[\Xsub{1}\ddd\Xsub{\rg}]
\gradn{\Ima{\Phg}}{\leq\kgtild}+
x_r\Ima{\Phg}.
\end{equation*}
Hence $\Ima{\Phg}=
\Rstar[\Xsub{1}\ddd\Xsub{\rg}]
\gradn{\Ima{\Phg}}{\leq\kgtild}$.
\item
When $p=Q$, we have $\gradn{\Ima{\Phg}}{\leq\kgtild}=0$,
so that $\Ima{\Phg}=0=N_Q$ by \eqref{cond25}. 
Next consider the case $\en{1}{p}{Q}$.
Since $\kgtild\leq-\kg_p$ by hypothesis and \eqref{cond25} holds, 
$\Ima{\Phg}$ satisfies the required conditions for 
the assertion of \Cor{\ref{proc13}} 
to be true. Hence 
\begin{equation*}
\Ima{\Phg}\sset
\Rstar[\Xsub{1}\ddd\Xsub{\rg}]
\gradn{\Ima{\Phg}}{\leq-\kg_p}\sset
N_p.
\end{equation*}
In consequence, the map $\pg$ restricted to $\dss{i=1}{p}\Rstar(\kg_i)$
is an injective homomorphism into $E$, which implies that
$\len{R}{E}=\infty$.
\item
Notice that $-\kg_1=\max\{\ d\ |\ \gradn{E/x_rE}{d}\neq0\ \}$.
Our assertion follows from the previous one.
\qed
\end{parenuma}
\killqed
\end{pf}

\section{Pseudo Weierstrass basis}
\lb{PWbasis}

Throughout this section we work over a \fgg\ module $E$ over $R$.
For a subset $E''$ of $E$,
let $(\seq{x}{j}{r})E'$ denote the set $x_jE'+\cdots+ x_rE'$ for
$\ee{1}{j}{r}$ and $(\seq{x}{j}{r})E':=\{0\}$ for $j=r+1$.

\begin{df}\lb{proc47}
Let $P=\{\ e^i_l\ |\ \ee{1}{i}{r+1},\ \ee{1}{l}{m_i}\ \}$ be a
set of homogeneous generators of $E$,
\begin{align}
&\supang{E}{i}: = \smm{l=1}{m_i} k[\seq{x}{i}{r}]e_l^i
\quad\fallt\quad \ee{1}{i}{r+1},\quad\andt
\lb{eq31}\\
&\supbr{E}{j}:=\smm{i=j}{r+1} \supang{E}{i}
\quad\fallt\quad \ee{1}{j}{r+1},
\lb{eq32}
\end{align}
where $ k[\seq{x}{i}{r}]=k$ for $i=r+1$ 
and $\smm{l=1}{m_i} ( \ \ )  = 0$ if $m_i = 0$.
We call $P$
a {\it pseudo \Wei\ basis} of $E$ with respect to $\seq{x}{1}{r}$, 
when it satisfies the following conditions.
\begin{parenumr}
\item\lb{cond18}
The module $k[\seq{x}{i}{r}] e^i_l$ is free over 
$k[\seq{x}{i}{r}]$ for each pair $i,l$.
\item\lb{cond19}
The sum \eqref{eq31} is direct over $k[\seq{x}{i}{r}]$.
\item\lb{cond20}
The sum \eqref{eq32} is direct over $k$.
\item\lb{cond21}
$E = \supbr{E}{1}$.
\item\lb{cond22}
For every triple $i',j,l$ satisfying
$1 \leq i' < j \leq r+1,\ \ee{1}{l}{m_j}$, we have
\begin{equation}
x_{i'}e^j_l  \in (\seq{x}{j}{r}) 
\lt(\smm{i=1}{j-1}\supang{E}{i}\rt)+\supbr{E}{j}.
\lb{eq33}
\end{equation}
\end{parenumr}
\end{df}

A \Wei\ basis of $E$ with respect 
to $\seq{x}{1}{r}$ defined in 
\cite[Section 2]{A5} is automatically a pseudo \Wei\ basis of 
$E$ with respect to $\seq{x}{1}{r}$. In the above definition,
the condition \eqref{cond22} is a loosened form of the
corresponding condition appearing in the definition of \Wei\ basis.

\begin{thm}\lb{proc48}
Suppose that the linear forms $\seq{x}{1}{r}$ form a reverse 
filter-regular $E$-sequence (see \cite[\Df{2.1}]{A5}).
Then there is a pseudo \Wei\ basis of $E$ with respect to
$\seq{x}{1}{r}$. In particular, there is a pseudo \Wei\ basis of 
$E$ with respect to $\seq{x}{1}{r}$ if $\Gg$ is chosen sufficiently 
generally.
\end{thm}

\begin{pf}
See \cite[\Thm{2.11}]{A5}.
\end{pf}

\begin{lem}\lb{proc52}
Let $P=\{\ e^i_l\ |\ \ee{1}{i}{r+1},\ \ee{1}{l}{m_i}\ \}$ be a
pseudo \Wei\ basis of $E$ with respect to $\seq{x}{1}{r}$.
Put $\Etild:=E/(\Hm{0}{E}+x_rE)$ and let 
$\etild^i_l:=e^i_l\ (\mod\ (\Hm{0}{E}+x_rE))$. Then 
$\{\ \etild^i_l\ |\ \ee{1}{i}{r},\ \ee{1}{l}{m_i}\ \}$
is a pseudo \Wei\ basis of the $R/(x_r)$-module 
$\Etild$ with respect to $\seq{x}{1}{r-1}$. 
\end{lem}

\begin{pf}
Left to the readers. See \cite[\Lem{2.4}, (3)]{A5} in which
the notation is a little bit different from the present one.
\end{pf}

\begin{lem}\lb{proc49}
Suppose that $\Gg$ is sufficiently general and 
let $P=\{\ e^i_l\ |\ \ee{1}{i}{r+1},\ \ee{1}{l}{m_i}\ \}$ be a
pseudo \Wei\ basis of $E$ with respect to $\seq{x}{1}{r}$.
Let further $(\nbar^1;\ldots;\nbar^{r+1})$ be the basic
sequence of $E$ (see \cite[\Df{2.13}]{A5}). 
Then $\nbar^i=(\deg(e^i_1)\ddd\deg(e^i_{m_i}))$
up to permutation for all $\ee{1}{i}{r+1}$.
\end{lem}

\begin{pf}
By \eqref{cond22} of \Df{\ref{proc47}}, $\supang{E}{r+1}=\supbr{E}{r+1}$
is a module over $R$. One sees therefore that
$\supang{E}{r+1}\sset\Hm{0}{E}\sset \un{t>0}0:_Ex_r^t\sset\supang{E}{r+1}$,
that is $\supang{E}{r+1}=\Hm{0}{E}$. This also holds for a \Wei\ basis of
$E$ with respect to $\seq{x}{1}{r}$ by \cite[(1) of \Lem{2.4}]{A5}.
Since $\{\ e^{r+1}_1\ddd e^{r+1}_{m_{r+1}}\ \}$ forms 
a basis of the vector space
$\Hm{0}{E}$ over $k$, one finds that  
$\nbar^{r+1}=(\deg(e^{r+1}_1)\ddd\deg(e^{r+1}_{m_{r+1}}))$
up to permutation. Consider the module 
$\Etild:=E/(\Hm{0}{E}+x_rE)$ over $R/(x_r)$. Then 
$\{\ \etild^i_l\ |\ \ee{1}{i}{r},\ \ee{1}{l}{m_i}\ \}$
is a pseudo \Wei\ basis of $\Etild$ with respect to $\seq{x}{1}{r-1}$ 
by \Lem{\ref{proc52}} and $(\nbar^1;\ldots;\nbar^{r})$ is the basic
sequence of $\Etild$ by \Thm{2.12} and \Df{2.13} of \cite{A5}. 
By induction we reach our assertion.
\end{pf}

\begin{lem}\lb{proc50}
Let $P=\{\ e^i_l\ |\ \ee{1}{i}{r+1},\ \ee{1}{l}{m_i}\ \}$ be a
pseudo \Wei\ basis of $E$ with respect to $\seq{x}{1}{r}$.
Then, 
\begin{equation}
R\supbr{E}{j+1}\sset\row{x}{j+1}{r}\lt(\smm{i=1}{j}\supang{E}{i}\rt)
+\supbr{E}{j+1}.
\lb{eq35}
\end{equation}
for every $\ee{1}{j}{r}$.
\end{lem}

\begin{pf}
Let us show our assertion by descending induction on $j$.
We see by \eqref{cond22} of \Df{\ref{proc47}} that
$\supbr{E}{r+1}$ is a submodule of $E$ over $R$.
This proves \eqref{eq35} for $j=r$. Suppose
that $\en{1}{j}{r}$ and that 
\begin{equation}
R\supbr{E}{j+2}\sset\row{x}{j+2}{r}
\lt(\smm{i=1}{j+1}\supang{E}{i}\rt)+\supbr{E}{j+2}.
\lb{eq34}
\end{equation}
Since $\supbr{E}{j+1}=\supang{E}{j+1}+\supbr{E}{j+2}$, we find
by \eqref{eq33} and \eqref{eq34} that
\begin{align*}
(\seq{x}{1}{r})\supang{E}{j+1}&\sset
(\seq{x}{j+1}{r})E+k[\seq{x}{j+1}{r}]\supbr{E}{j+1}\\
&\sset(\seq{x}{j+1}{r})E+\supang{E}{j+1}+R\supbr{E}{j+2}\\
&\sset(\seq{x}{j+1}{r})E+\supbr{E}{j+1},
\end{align*}
so that 
\begin{equation*}
(\seq{x}{1}{r})\supbr{E}{j+1}
\sset(\seq{x}{j+1}{r})E+\supbr{E}{j+1}
\end{equation*}
by \eqref{eq34}.
Repeated use of this inclusion yields 
\begin{equation*}
R\supbr{E}{j+1}
\sset(\seq{x}{j+1}{r})E+\supbr{E}{j+1}.
\end{equation*}
Since $E=\supbr{E}{1}$, we see that \eqref{eq35} holds for $j$
again by \eqref{eq34}. 
\end{pf}

\begin{lem}\lb{proc51}
Let $P=\{\ e^i_l\ |\ \ee{1}{i}{r+1},\ \ee{1}{l}{m_i}\ \}$ be a
pseudo \Wei\ basis of $E$ with respect to $\seq{x}{1}{r}$ and
$i_0,l_0$ be integers with $\ee{1}{i_0}{r},\ \ee{1}{l_0}{m_{i_0}}$.
Let further $\hg$ be a homogeneous element of $R\supbr{E}{i_0+1}$ of degree
$\deg(e^{i_0}_{l_0})$.
For each pair of integers $i,l$ with $\ee{1}{i}{r},\ \ee{1}{l}{m_{i}}$,
put $e'{}^i_l:=e^{i_0}_{l_0}+\hg$ if $(i,l)=(i_0,l_0)$ and 
$e'{}^i_l:=e{}^i_l$ otherwise. Then 
$P':=\{\ e'{}^i_l\ |\ \ee{1}{i}{r+1},\ \ee{1}{l}{m_i}\ \}$ is also
a pseudo \Wei\ basis of $E$ with respect to $\seq{x}{1}{r}$.
\end{lem}

\begin{pf}
Put 
\begin{equation*}
\supangp{E}{i}: = \smm{l=1}{m_i} k[\seq{x}{i}{r}]e'{}_l^i
\quad\fort\ \ee{1}{i}{r+1}\ \andt\quad
\supbrp{E}{j}:=\smm{i=j}{r+1} \supangp{E}{i}.
\end{equation*}
\begin{claim}[Claim 1]
$\supbrp{E}{1}=E$.
\end{claim}
\begin{pfofclaim}[Proof of Claim 1]
This claim can be proved by induction on $r-i_0$. 
Consider first the case $i_0=r$. Since 
$\supbrp{E}{i_0+1}=\supbr{E}{r+1}=\Hm{0}{E}$ is a module over $R$,
$R\hg\in\supbrp{E}{i_0+1}$. Hence 
$E=\supbr{E}{1}\sset\supbrp{E}{1}+R\hg
\sset\supbrp{E}{1}+\supbrp{E}{i_0+1}=\supbrp{E}{1}\sset E$.
Now suppose that $r-i_0>0$ and that our assertion is true for 
all smaller values of $r-i_0$.
Put $\Etild:=E/(\Hm{0}{E}+x_rE)$ and let $\hgtild$,
$\etild^i_l$ and $\etild{}'{}^i_l$ denote the elements of
$\Etild$ obtained from $\hg$, $e^i_l$ and $e'{}^i_l$ respectively modulo
$\Hm{0}{E}+x_rE$. Put further
\begin{equation*}
\supangp{\Etild}{i}: = \smm{l=1}{m_i} k[\seq{x}{i}{r-1}]\etild{}'{}_l^i
\quad\fort\ \ee{1}{i}{r}\ \andt\quad
\supbrp{\Etild}{j}:=\smm{i=j}{r} \supangp{\Etild}{i}.
\end{equation*}
Then the set $\{\ \etild^i_l\ |\ \ee{1}{i}{r},\ \ee{1}{l}{m_i}\ \}$
is a pseudo \Wei\ basis of the module $\Etild$ over $R/(x_r)$ and
$\hgtild\in(R/(x_r))\supbrp{\Etild}{i_0+1}$. Since
$(r-1)-i_0<r-i_0$, we have $\Etild=\supbrp{\Etild}{1}$ by our
induction hypothesis. In other words, $E=\supbrp{E}{1}+x_rE$.
Since $x_r\supbrp{E}{1}\sset\supbrp{E}{1}$, plugging the above
equality into $E$ on the right hand side successively, we see that
\begin{align*}
E&=\supbrp{E}{1}+x_r(\supbrp{E}{1}+x_rE)\\
&=\supbrp{E}{1}+x_r^2E=\supbrp{E}{1}+x_r^\ng E
\end{align*}
for all $\ng>0$. Thus $\supbrp{E}{1}=E$.
\qedclaim
\end{pfofclaim}
\begin{claim}[Claim 2]
For all $\ee{1}{j}{r}$, we have
\begin{equation}
R\supbrp{E}{j+1}\sset\row{x}{j+1}{r}\lt(\smm{i=1}{j}\supangp{E}{i}\rt)
+\supbrp{E}{j+1}.
\lb{eq36}
\end{equation}
\end{claim}
\begin{pfofclaim}[Proof of Claim 2]
Let the notation be as in the proof of Claim 1. We again 
show our assertion by induction on $r-i_0$. In the case $r=i_0$,
$\supangp{E}{i}=\supang{E}{i}$ and $\supbrp{E}{i}=\supbr{E}{i}$
for $\en{1}{i}{r}$, $\supang{E}{r}\sset\supangp{E}{r}+R\hg$,
$\supangp{E}{r}\sset\supang{E}{r}+R\hg$, and
$R\hg\in\supbrp{E}{i_0+1}=\supbr{E}{r+1}=\Hm{0}{E}$. 
It follows therefore that
\begin{align*}
R\supbrp{E}{j+1}
&\sset R\supbr{E}{j+1}+R\hg\\
&\sset\row{x}{j+1}{r}\lt(\smm{i=1}{j}\supang{E}{i}\rt)
+\supbr{E}{j+1}+\supbrp{E}{r+1}\\
&\sset\row{x}{j+1}{r}\lt(\smm{i=1}{j}\supangp{E}{i}\rt)
+\supbrp{E}{j+1}+R\hg+\supbrp{E}{r+1}\\
&\sset\row{x}{j+1}{r}\lt(\smm{i=1}{j}\supangp{E}{i}\rt)
+\supbrp{E}{j+1}
\end{align*}
for all $\ee{1}{j}{r}$.
Suppose that $r-i_0>0$ and that \eqref{eq36} holds for
all smaller values of $r-i_0$. Then, since $E=\supbrp{E}{1}$ and
\begin{equation*}
(R/(x_r))\supbrp{\Etild}{j+1}
\sset\row{x}{j+1}{r-1}\lt(\smm{i=1}{j}\supangp{\Etild}{i}\rt)
+\supbrp{\Etild}{j+1}
\end{equation*}
for all $\ee{1}{j}{r-1}$ by the induction hypothesis, 
\begin{align*}
R\supbrp{E}{j+1}
&\sset\row{x}{j+1}{r-1}\lt(\smm{i=1}{j}\supangp{E}{i}\rt)
+\supbrp{E}{j+1}+\Hm{0}{E}+x_rE\\
&\sset\row{x}{j+1}{r}\lt(\smm{i=1}{j}\supangp{E}{i}\rt)
+\supbrp{E}{j+1}
\end{align*}
for all $\en{1}{j}{r}$. Besides, the inclusion \eqref{eq36} is
clear for $j=r$.
\qedclaim
\end{pfofclaim}
\begin{claim}[Claim 3]
$\smm{i=1}{r+1}\smm{l=1}{m_i}g^i_le'{}^i_l=0$ with 
$g^i_l\in k[\seq{x}{i}{r}]$ for all $i,l$ if and 
only if $g^i_l=0$ for all $i,l$.
\end{claim}
\begin{pfofclaim}[Proof of Claim 3]
By \eqref{cond18} -- \eqref{cond21} of \Df{\ref{proc47}}, 
\begin{equation*}
\dimk{\gradn{E}{\ng}}
=\smm{i=1}{r+1}\smm{l=1}{m_i}\dimk{\gradn{k[\seq{x}{i}{r}]}{\ng-\deg(e^i_l)}}
\quad\fallt\quad\ng\in\Zbf.
\end{equation*}
On the other hand, 
\begin{equation*}
\dimk{\gradn{\supbrp{E}{1}}{\ng}}\leq
\smm{i=1}{r+1}\smm{l=1}{m_i}\dimk{\gradn{k[\seq{x}{i}{r}]}{\ng-\deg(e'{}^i_l)}}
\quad\fallt\quad\ng\in\Zbf.
\end{equation*}
Moreover $\deg(e'{}^i_l)=\deg(e^i_l)$ for all $i,l$ by the
definition of $e'{}^i_l$. 
If there were a nontrivial relation 
$\smm{i=1}{r+1}\smm{l=1}{m_i}g^i_le'{}^i_l=0$
with $g^i_l\in k[\seq{x}{i}{r}]$ for all $i,l$, then 
$\dimk{\gradn{\supbrp{E}{1}}{\ng}}<\dimk{\gradn{E}{\ng}}$ for some
$\ng\in\Zbf$.
This contradicts Claim 1. Hence our assertion holds.
\qedclaim
\end{pfofclaim}
The conditions \eqref{cond18} -- \eqref{cond22} of \Df{\ref{proc47}}
are all satisfied by $P'$ by Claims 1 -- 3.
\end{pf}

\section{Gorenstein graded rings of dimension zero}
\lb{gor}

\par 
In this section, we show how we can apply the results of
section \ref{mat} to Gorenstein graded rings.

\par
Let $I$ be a homogeneous ideal in $R$ 
such that the graded ring $A:=R/I$ is of dimension zero and
let $B_R(I)=(a;n^2_1\ddd n^2_a;\cdots;n^r_1\ddd n^r_{m_r})$ be
its basic sequence.
Notice that $m_{r+1}=0$ and $\nbar^{r+1}=\eset$ since $\Hm{0}{I}=0$.
Assume that $r\geq3$ and that $\Gg$ is chosen sufficiently generally.
Put  $\rg:=r-1$ and $\Rstar:=k[x_r]$.
Let $P=\{\ f^i_l\ |\ \ee{1}{i}{r},\ \ee{1}{l}{m_i}\ \}$ be 
a pseudo \Wei\ basis of
$I$ with respect to $\seq{x}{1}{r}$,
$\supang{I}{i}:=\dss{l=1}{m_i}k[\seq{x}{i}{r}]f^i_l$ $(\ee{1}{i}{r})$ 
(see the previous section), and $I':=\smm{i=1}{\rg}\supang{I}{i}$. 
Note that $I=\smm{i=1}{r}\supang{I}{i}=I'+\supang{I}{r}$ and that these
sums are direct over $\Rstar$.
For the sake of convenience, we 
put $\Rbar:=R/(x_r)=k[\seq{x}{1}{r-1}]$, 
$\nfrak:=(\seqsub{x}{1}{r-1})\sset k[\seq{x}{1}{r-1}]$, 
$\Rbarstar:=\Rstar/x_r\Rstar$, $\Abar:=A/x_rA$, 
$\Ibar:=I+(x_r)/(x_r)\sset\Rbar$ and $\Itild:=I/x_rI$. 
For an element $f\in I$, put $\ftild:=f\ (\mod\ x_rI)\in I/x_rI$
and $\fbar:=f\ (\mod\ (x_r))\in \Ibar$.

\begin{lem}\lb{proc40}
With the notation above, $R/I'$ is a \fgg\ free module over $\Rstar$
of rank $m_r$, $\supang{I}{r}\sset(x_r)$, $I'\cap(x_r)=x_rI'$, 
and the ring $A$ has a minimal 
free resolution over $\Rstar$ of the form
\begin{equation}
0\lra\dss{i=1}{m_r}\Rstar(-n^r_i)\xra{\ \Phg\ }
\dss{i=1}{m_r}\Rstar(\kg_i)
\xra{\ \pg\ }A\lra 0
\lb{eq24}
\end{equation}
with homogeneous homomorphisms $\Phg$ and $\pg$ of degree zero.
Here, $\dss{i=1}{m_r}\Rstar(\kg_i)\cong R/I'$ as $\Rstar$-modules 
and the integers $\seq{\kg}{1}{m_r}$
satisfy $\kg_l\leq \kg_{l+1}\leq\kg_l+1$ for
all $\en{1}{l}{m_r}$ and $\kg_{m_r}=0$.
\end{lem}

\begin{pf}
The exact sequence 
\begin{equation}
0\lra I\lra R\lra A\lra 0
\lb{eq37} 
\end{equation}
induces an exact sequence
\begin{equation}
0\lra\Tor{1}{\Rstar}{A}{\Rbarstar}\xra{\ \eg'''\ } I/x_rI
\xra{\ \eg\ } R/(x_r)\lra A/x_rA\lra0.
\lb{eq25}
\end{equation}
By \Lem{\ref{proc52}}, the set 
$\{\ \ftild^i_l\ |\ \ee{1}{i}{r},\ \ee{1}{l}{m_i}\ \}$ is
a pseudo \Wei\ basis of $I/x_rI$ with respect to $\seq{x}{1}{r-1}$ 
and $\supang{I}{r}/x_r\supang{I}{r}
\cong\dss{l=1}{m_r}k\ftild^r_l=\supang{I}{r}+x_rI/x_rI
=\Hsupb{0}{\nfrak}{I/x_rI}
=\Ker{\eg}\cong\Tor{1}{\Rstar}{A}{\Rbarstar}$. 
In particular, $f^r_l\in(x_r)$ for all
$\ee{1}{l}{m_r}$. Hence $\supang{I}{r}\sset(x_r)$. 
Moreover, $I'\cap(x_r)\sset\supang{I}{r}+x_rI=\supang{I}{r}+x_rI'$ since 
$I'\cap(x_r)\ (\mod\ x_rI)\sset\Ker{\eg}$.
Hence $I'\cap(x_r)=x_rI'$.
On the other hand, the exact sequence 
\begin{equation}
0\lra\supang{I}{r}=I/I'\lra R/I'\lra A\lra 0
\lb{eq26}
\end{equation}
yields another exact sequence
\begin{align}
0\lra\Tor{1}{\Rstar}{R/I'}{\Rbarstar}
&\lra\Tor{1}{\Rstar}{A}{\Rbarstar}
\xra{\ \eg''\ }\supang{I}{r}/x_r\supang{I}{r}
\notag\\
&\lra R/I'+(x_r)\xra{\ \eg'\ } A/x_rA\lra0.
\lb{eq27}
\end{align}
One has the natural morphism over $\Rstar$ from the complex \eqref{eq37} to
\eqref{eq26}, which induces a morphism from \eqref{eq25} to
\eqref{eq27}.
The map $\eg'''$ followed by the natural surjection
from $I/x_rI$ to $I/x_rI+I'=\supang{I}{r}/x_r\supang{I}{r}$, therefore,
coincides with $\eg''$.
Since $\Ima{\eg'''}=\Ker{\eg}=\supang{I}{r}+x_rI/x_rI$, 
it follows that $\Ker{\eg''}=\Ker{\eg'''}=0$.
On the other hand,  the map $\eg'$ is an isomorphism 
since $\supang{I}{r}\sset(x_r)$ as we have already seen, so that
$\eg''$ is surjective.
In consequence, one finds that $\Tor{1}{\Rstar}{R/I'}{\Rbarstar}=0$. 
Hence $R/I'$ is
a \fgg\ free module over $\Rstar$. Since $\len{\Rstar}{A}<\infty$,
the rank of $R/I'$ over $\Rstar$ coincides with that of
$\supang{I}{r}$. The exact sequence \eqref{eq26} gives the
desired minimal free resolution \eqref{eq24}.
With no loss of generality, we can choose $\seq{\kg}{1}{m_r}$ so 
that they form a nondecreasing sequence. 
Suppose $\kg_{l_0+1}>\kg_{l_0}+1$ for some $\en{1}{l_0}{m_r}$. 
Then $\gradn{R}{1}\gradn{A}{\leq -\kg_{l_0+1}}
\sset\gradn{A}{\leq -\kg_{l_0+1}+1}
\sset\Rstar\gradn{A}{\leq -\kg_{l_0+1}}$ by \eqref{eq24}. Hence 
$R\gradn{A}{\leq -\kg_{l_0+1}}
\sset\Rstar\gradn{A}{\leq -\kg_{l_0+1}}\neq A$. 
Since $1\in\gradn{A}{0}$, the largest $\kg_{m_r}$ must be $0$.
This means that $1\in\gradn{A}{\leq -\kg_{l_0+1}}$, therefore
$R\gradn{A}{\leq -\kg_{l_0+1}}=A$, which is a contradiction.
Thus $\kg_l\leq \kg_{l+1}\leq\kg_l+1$ for
all $\en{1}{l}{m_r}$.
\end{pf}

\begin{lem}\lb{proc44}
With the notation above,
\begin{equation*}
-\kg_1=n^{r-1}_{m_{r-1}}-1=\max\{\ d\ |\ \gradn{A/x_rA}{d}\neq0\ \}.
\end{equation*}
\end{lem}

\begin{pf}
It follows from the exact sequence \eqref{eq25} that 
$\Ibar=\Ima{\eg}\cong \Itild/\Hn{0}{\Itild}$.
Observe that the basic sequences of $I$ and $\Itild$ are the same and that
the basic sequence of $\Ibar$ is 
$(a;n^2_1\ddd n^2_a;\cdots;n^{r-1}_1\ddd n^{r-1}_{m_{r-1}})$
(cf. \cite[\Lem{2.4}]{A5}). By \Lem{\ref{proc40}} applied to
$\Abar$, there is a minimal free resolution of $\Abar$
over $\Rstar$ of the form
\begin{equation*}
0\lra\dss{i=1}{m_{r-1}}\Rstar(-n^{r-1}_i)\lra
\dss{i=1}{m_{r-1}}\Rstar(\kgtild_i)
\lra\Abar\lra 0,
\end{equation*}
so that $\min\{\ d\ |\ \gradn{\Ext{1}{\Rstar}{\Abar}{\Rstar}}{d}\neq0\ \}
=-n^{r-1}_{m_{r-1}}$. 
In addition,  
$\Ext{1}{\Rstar}{\Abar}{\Rstar}\cong\Hom{k}{\Abar}{k}(1)$
by duality. Hence 
$\min\{\ d\ |\ \gradn{\Hom{k}{\Abar}{k}}{d}\neq0\ \}
=-n^{r-1}_{m_{r-1}}+1$. On the other hand,
$\max\{\ d\ |\ \gradn{\Abar}{d}\neq0\ \}=-\kg_1$ by \eqref{eq24}.
Thus $-\kg_1=-(-n^{r-1}_{m_{r-1}}+1)$, which proves our assertion.
\end{pf}

\begin{lem}\lb{proc53}
With the notation of \Lem{\ref{proc40}}, 
let $\{\seq{e}{1}{m_r}\}$ be
the canonical basis of $\dss{i=1}{m_r}\Rstar(\kg_i)$, and
for each $\ee{1}{i}{m_r}$,
let $\gcheck_i$
denote a homogeneous element of $R$ satisfying $e_i=\gcheck_i\ (\mod\ I')$
via the isomorphism $\dss{i=1}{m_r}\Rstar(\kg_i)\cong R/I'$. 
Denoting by $\phgsub{ij}$ the $(i,j)$-component of
$\Phg$, put $e^r_l:=\smm{i=1}{m_r}\phgsub{il}\gcheck_i\in I$.
Then $\Pcheck:=\{\ f^i_l\ |\ \ee{1}{i}{r-1},\ \ee{1}{l}{m_i}\ \}
\cup\{ e^r_l\ |\ \ee{1}{l}{m_r}\ \}$ is a pseudo \Wei\ basis of $I$
with respect to $\seq{x}{1}{r}$.
\end{lem}

\begin{pf}
By its construction, $\Pcheck$ satisfies the conditions 
\eqref{cond18} -- \eqref{cond21} of \Df{\ref{proc47}}.
On the other hand $e^r_l\in\supang{I}{r}+x_rI=x_rI'+\supbr{I}{r}$ 
by \eqref{eq25},
since $e^r_l$ lies in $I\cap(x_r)$, i.e. $e^r_l\ (\mod\ x_rI) \in \Ker{\eg}$.
Moreover, since $I'\cap(x_r)=x_rI'$ and 
since both $\{e^r_1\ddd e^r_{m_r}\}$ and $\{f^r_1\ddd f^r_{m_r}\}$
modulo $I'$ form free bases of $I/I'$ over $\Rstar$, it follows that
$f^r_l\in\smm{i=1}{m_r}\Rstar e^r_i+x_rI'$. With the use of these 
observations the verification of the condition \eqref{cond22} of 
\Df{\ref{proc47}} can now easily be carried out for $\Pcheck$.
\end{pf}

\par
The next lemma will be used afterwards in Section \ref{Ibar}.
 
\begin{lem}\lb{proc54}
With the notation of \Lem{\ref{proc53}}, assume that 
$f^r_l=\smm{i=1}{m_r}\phgsub{il}\gcheck_i$ for all
$\ee{1}{l}{m_r}$ from the first, namely $P=\Pcheck$.
Put $P_1:=P-\{\ f^r_l\ |\ \ee{1}{l}{m_r}\ \}$.
Suppose further that there are $\seq{f}{1}{s}\in P_1$ satisfying the following
two conditions.
\begin{parenumr}
\item\lb{cond23}
The polynomials $\seq{\fbar}{1}{s}$ generate $\Ibar$
minimally over $\Rbar$. 
\item\lb{cond24}
$P_1-\{\seq{f}{1}{s}\}\sset\smm{i=1}{s}Rf_i$.
\end{parenumr}
Let $\Xsub{i}$ $(\ee{1}{i}{\rg})$ be the matrices 
and $\thgsub{ijl}$ $(\ee{1}{l}{m_r})$ the columns of 
$[\Xsub{i},\Xsub{j}]$ mentioned in Section \ref{mat}
associated with \eqref{eq24}.
If $\ee{1}{l}{m_r}$ and 
\begin{equation}
\dimk{\gradplain{\Hn{0}{\Itild}/\nfrak\Hn{0}{\Itild}}{-\kg_l+2}}
+\dimk{\gradplain{\Ibar/\nfrak\Ibar}{-\kg_l+2}}
=\dimk{\gradplain{\Itild/\nfrak\Itild}{-\kg_l+2}},
\lb{eq38}
\end{equation}
then $\thgsub{ijl}\in
\smm{i=1}{\rg}\Xsub{i}\gradn{\Ima{\Phg}}{-\kg_l+1}
+x_r\gradn{\Ima{\Phg}}{-\kg_l+1}$.
\end{lem}

\begin{pf}
Let $\seq{f}{s+1}{s+t}$ be elements of 
$\{\ f^r_l\ |\ \ee{1}{l}{m_r}\ \}$ such that $\seq{\ftild}{s+1}{s+t}$
minimally generate $\Hn{0}{\Itild}$ over $\Rbar$.
Since there is an exact sequence 
\begin{equation*}
0\lra\Hn{0}{\Itild}\lra\Itild\lra\Ibar\lra0
\end{equation*}
by \eqref{eq25} and since the condition \eqref{cond23} holds, 
it follows that $\Itild$ is generated by $\seq{\ftild}{1}{s+t}$
over $\Rbar$. Hence $I\sset(\seq{f}{1}{s+t})$.  
Moreover 
\begin{equation*}
\supang{I}{r}\sset\smm{j=1}{t}R f_{s+j}+x_r\supang{I}{r}+x_rI'\sset
\smm{j=1}{t}R f_{s+j}+x_r\smm{i=1}{s}Rf_i,
\end{equation*}
since $I'\sset(\seq{f}{1}{s})$ by the condition \eqref{cond24}.
Now let $\phi_{i,\mg\ng}$ denote the $(\mg,\ng)$-component of $\Xsub{i}$.
We have $x_j\gcheck_l\in\smm{\ng=1}{m_r}\gcheck_\ng\phi_{j,\ng l}+I'$ and
$x_ix_j\gcheck_l
\in\smm{\mg=1}{m_r}\smm{\ng=1}{m_r}
\gcheck_\mg\phi_{i,\mg\ng}\phi_{j,\ng l}+RI'$.
Likewise, $x_jx_i\gcheck_l
\in\smm{\mg=1}{m_r}\smm{\ng=1}{m_r}
\gcheck_\mg\phi_{j,\mg \ng}\phi_{i,\ng l}+RI'$.
After subtraction, 
\begin{equation*}
f:=(\seq{\gcheck}{1}{m_r})\thgsub{ijl}
=\smm{\mg=1}{m_r}
\gcheck_\mg\smm{\ng=1}{m_r}(\phi_{i,\mg \ng}\phi_{j,\ng l}-
\phi_{j,\mg \ng}\phi_{i,\ng l})\in RI'=(\seq{f}{1}{s}).
\end{equation*}
On the other hand, $\thgsub{ijl}\in\Ima{\Phg}$,
so that $f\in\supang{I}{r}$. Thus we can write
$f=\smm{i=1}{s}h_if_i+\smm{j=1}{t}h'_jf_{j+s}=\smm{i=1}{s}h''_if_i$
with homogeneous $h_i\in x_rR$ and $h'_j,\ h''_i\in R$. 
Notice that $\deg{f}=-\kg_l+2$.
Our assumption \eqref{eq38} means that the elements of degree $-\kg_l+2$
taken from $\seq{\ftild}{1}{s+t}$ are
linearly independent over $k$ when considered modulo $\nfrak\Itild$. 
Since $\smm{i=1}{s}(h_i-h''_i)f_i+\smm{j=1}{t}h'_jf_{s+j}=0$,
we must have $h_i-h''_i$ and $h'_j$ lie in $\mfrak$ for all $i,j$.
Moreover, 
$R\supang{I}{r}=R\supbr{I}{r}\sset x_rI'+\supbr{I}{r}
=x_rI'+\supang{I}{r}$
by \Lem{\ref{proc50}} and $RI'\sset I'+\supang{I}{r}$.
Hence 
\begin{align*}
f\in x_rRI'&+(\seq{x}{1}{r})R\supang{I}{r}
\sset x_rI'+(\seq{x}{1}{r})x_rI'+(\seq{x}{1}{r})\supang{I}{r}\\
&\sset(\seq{x}{1}{r})\supang{I}{r}+x_rI'\\
&\sset(\seq{\gcheck}{1}{m_r})\lt(\smm{i=1}{\rg}
\Xsub{i}\Ima{\Phg}+x_r\Ima{\Phg}\rt)+I'.
\end{align*}
This implies our assertion.
\end{pf}

\begin{df}\lb{proc62}
A graded ring $A$ has the
weak Lefschetz property (WLP) if there is an element $\ell\in\gradn{A}{1}$
such that 
\begin{equation*}
\tm\ell:\gradn{A}{i}\lra\gradn{A}{i+1}
\end{equation*}
has maximal rank for every $i\in\Zbf$. $A$ has the strong Lefschetz property
(SLP) if there is an element $\ell\in\gradn{A}{1}$
such that 
\begin{equation*}
\tm\ell^d:\gradn{A}{i}\lra\gradn{A}{i+d}
\end{equation*}
has maximal rank for every $i\in\Zbf$ and $d\in\Nbf$.
\end{df}

We are concerned only with WLP in this paper.
In the case $\ch{k}=0$, another proof is available for the
lemma below in \cite[\Prop{3.4}]{ACP}.

\begin{lem}\lb{proc5}
The graded ring $A$ has the WLP if and only if 
$n^r_1\geq n^{r-1}_{m_{r-1}}$.
\end{lem}

\begin{pf}
Since $\max\{\ -\kg_l\ |\ \ee{1}{l}{m_r}\ \}=-\kg_1$ and
$\min\{\ n^r_l\ |\ \ee{1}{l}{m_r}\ \}=n^r_1$, the map 
$\tm x_r:\gradn{A}{i}\lra\gradn{A}{i+1}$ 
is surjective for all $i\geq-\kg_1$ and injective for all $i < n^r_1-1$
by \Lem{\ref{proc40}}.
If $n^r_1\geq r^{r-1}_{m_{r-1}}$, then $-\kg_1=n^{r-1}_{m_{r-1}}-1
\leq n^r_1-1$, it has therefore maximal rank for all $i\in\Zbf$, 
that is, $A$ has the WLP.
On the contrary, suppose that $n^r_1<r^{r-1}_{m_{r-1}}$. 
Then $-\kg_1\geq n^r_1$.
Since $\kg_l\leq \kg_{l+1}\leq\kg_l+1$ for all $\en{1}{l}{m_r}$ and
$\kg_{m_r}=0$, 
there is an $l_1$ such that $-\kg_{l_1}=n^r_1$. In this case
the map $\tm x_r:\gradn{A}{n^r_1-1}\lra\gradn{A}{n^r_1}$ can
neither be surjective nor injective. Since the linear form $x_r$ is 
sufficiently general, this means that $A$ does not have the WLP.
\end{pf}

\par
Put $\socd{A}:=\max\{\ d\ |\ \gradn{A}{d}\neq0\ \}$.

\begin{lem}\lb{proc37}
With the notation above, suppose further that $A$ is \Gor.
Then, $\socd{A}=n^r_{m_r}-1$ and 
there is a minimal free resolution of 
$A$ over $\Rstar$ of the form
\begin{equation}
0\lra\dss{i=1}{m_r}\Rstar(-n^r_i)\xra{\ \Phg\ }
\dss{i=1}{m_r}\Rstar(-n^r_{m_r}+n^r_i)
\xra{\ \pg\ }A\lra 0
\lb{eq9}
\end{equation}
with homogeneous homomorphisms $\Phg$ and $\pg$ of degree zero.
\end{lem}

\begin{pf}
By \Lem{\ref{proc40}}, one sees that
$\Ext{1}{\Rstar}{A}{\Rstar}$ is a factor module of 
$\dss{l=1}{m_r}\Rstar(n^r_l)$
by its submodule contained in $x_r\dss{l=1}{m_r}\Rstar(n^r_l)$.
This implies that 
\begin{equation*}
\min\{\ d\ |\ \gradn{\Ext{1}{\Rstar}{A}{\Rstar}}{d}\neq0\ \}
=-n^r_{m_r}.
\end{equation*}
On the other hand, 
\begin{equation*}
\min\{\ d\ |\ \gradn{\Hom{k}{A}{k}}{d}\neq0\ \}=-\socd{A}
\end{equation*}
and $\Ext{r}{R}{A}{R}(-r)\cong\Hom{k}{A}{k}
\cong\Ext{1}{\Rstar}{A}{\Rstar}(-1)$.
Hence $\socd{A}=n^r_{m_r}-1$.
Since $\Ext{1}{\Rstar}{A}{\Rstar}\cong\Hom{k}{A}{k}(1)
\cong A(\socd{A}+1)$, it follows from the dual of \eqref{eq24} over
$\Rstar$ that 
$(\kg_1\ddd\kg_{m_r})=(-\socd{A}-1+n^r_1\ddd-\socd{A}-1+n^r_{m_r})
=(-n^r_{m_r}+n^r_1\ddd -n^r_{m_r}+n^r_{m_r})$.
Thus one obtains the exact sequence \eqref{eq9}.
\end{pf}

\begin{rem}\lb{proc38}
In the case $A$ is \Gor, one has $0:_A \mfrak\cong k(-\socd{A})$ and
\begin{equation}
A(n^r_{m_r}-1)\cong\Ext{r}{R}{A}{R}(-r)\cong\Ext{1}{\Rstar}{A}{\Rstar}(-1)
\cong\Hom{k}{A}{k},
\lb{eq49}
\end{equation}
as mentioned in the above proof. Moreover,
\begin{equation}
(0:_Ax_r)(\socd{A})=\Hom{k}{\Abar}{k}
=\Ext{r-1}{\Rbar}{\Abar}{\Rbar}(-r+1).
\lb{eq50}
\end{equation}
\end{rem}

\begin{cor}\lb{proc39}
If $A$ is \Gor,
\begin{equation*}
n^r_{m_r}-n^r_1=n^{r-1}_{m_{r-1}}-1=\max\{\ d\ |\ \gradn{A/x_rA}{d}\neq0\ \}.
\end{equation*}
\end{cor}

\begin{pf}
Clear by \Lems{\ref{proc44} and \ref{proc37}}.
\end{pf}

\par
Let $\seq{\kg}{1}{m_r}$, $\Phg$ and $\pg$ be as in \Lem{\ref{proc40}},
and $\Xsub{j}$ $(\ee{1}{j}{\rg})$ be the matrices 
as stated in Section \ref{mat}. We identify $\Phg$ with the matrix
representing itself.
Put $L:=\dss{i=1}{m_r}\Rstar(\kg_i)$ and
$F:=\dss{i=1}{m_r}\Rstar(-n^r_i)$.
For each $\ee{1}{i}{\rg}$, let $\Ysub{i}$ be a matrix in $\Rstar$
which gives a homogeneous homomorphism 
$\Ysub{i}:F(-1)\lra F$
satisfying
$\Xsub{i}\Phg=\Phg\Ysub{i}$.
Since $\pg\Xsub{i}\Phg=x_i\pg\Phg=0$,
such a $\Ysub{i}$ exists indeed for each $i$ and
is determined uniquely by $\Xsub{i}$ and $\Phg$.
They make the following diagram commutative :
\begin{equation}
\begin{CD}
F(-1) @>{\ \Phg\ }>>
L(-1) @>{\ \pg\ }>>A(-1)\\
@V \Ysub{i} VV @VV \Xsub{i} V @VV \tm x_i V\\
F @>{\ \Phg\ }>>
L @>{\ \pg\ }>>A
\end{CD}\qquad.
\lb{eq14}
\end{equation}

\par
For a module $M$ over $\Rstar$, we will denote by
$M^\vee$ the dual module $\Hom{\Rstar}{M}{\Rstar}$.

\begin{rem}\lb{proc7}
Suppose that $A$ is \Gor.
\begin{parenuma}
\item\lb{cond36}
Since the dual of \eqref{eq9} gives a minimal free resolution
of the $\Rstar$-module 
$A(n^r_{m_r})\cong\Ext{1}{\Rstar}{A}{\Rstar}$ over $\Rstar$, 
we may assume that $\Ima{\Phg}
=\Imasup{\Rstar}{\tpose{\Phg}}$.
Moreover, since the isomorphism $A(n^r_{m_r})\cong\Ext{1}{\Rstar}{A}{\Rstar}$
is also an isomorphism over $R$, the multiplication by 
$\tpose{\Ysub{i}}$ of an element of $F^\vee$ represents
the multiplication by $x_i$ of the corresponding element of $A$. 
We can therefore choose $\Xsub{i}$, $\Ysub{i}$ 
$(\ee{1}{i}{\rg})$ and $\Phg$ so that every column of 
$\tpose{\Ysub{i}}-\Xsub{i}$
lies in $\Ima{\Phg}$.
In particular, $\tpose{\Ysub{i}}\equiv\Xsub{i}\ 
(\mod\ x_r\Rstar)$.
\item\lb{cond37}
Without changing $F$ and $L$, one can normalize $\Phg$ to a 
symmetric matrix such that there is only one nonzero component in each 
row and column. But when that is done, one may only be able to say
that every column of $C^{-1}\tpose{\Ysub{i}}C-\Xsub{i}$
lies in $\Ima{\Phg}$ with an invertible matrix 
$C$. Moreover, $A$ has the strong Lefschetz property if and only if
a normal form of $\Phg$ mentioned above is a diagonal matrix.
\end{parenuma}
\end{rem}

\par
The canonical basis of $F$ (resp. $L$)
over $\Rstar$ will be denoted by 
$\{e^F_1\ddd e^F_{m_r}\}$ (resp. $\{e^L_1\ddd e^L_{m_r}\}$).
Further, the canonical basis of $L^\vee$ will be denoted by
$\{e^{L\ast}_1\ddd e^{L\ast}_{m_r}\}$ in a standard way.
Given an integer $u$ with $\ee{1}{u}{m_r}$, let 
\begin{align*}
&M_{1,u}:=\smm{l=1}{u}
\Rstar[\Xsub{1}\ddd\Xsub{\rg}]\Phg e^F_l,\\
&M_{2,u}:=\smm{l=1}{u}
\Rstar[\tpose{\Xsub{1}}\ddd\tpose{\Xsub{\rg}}] e^{L\ast}_l,\quad
M_{3,u}:=\smm{l=1}{u}
\Rstar[\Ysub{1}\ddd\Ysub{\rg}] e^F_l.
\end{align*}
Then,
$M_{1,u}$ is a submodule of $L$ over the $\Rstar$-algebra
$\Rstar[\Xsub{1}\ddd\Xsub{\rg}]$, 
$M_{2,u}$ a submodule of $L^\vee$ over
$\Rstar[\tpose{\Xsub{1}}\ddd\tpose{\Xsub{\rg}}]$,
and $M_{3,u}$ a submodule of $F$ over
$\Rstar[\Ysub{1}\ddd\Ysub{\rg}]$.
Notice that $\Rstar[\Xsub{1}\ddd\Xsub{\rg}]\Phg e^F_l
=\Phg\Rstar[\Ysub{1}\ddd\Ysub{\rg}] e^F_l$
for all $\ee{1}{l}{m_r}$ by the commutativity of \eqref{eq14}.

\begin{lem}\lb{proc12}
With the notation above, we have
$M_{1,u}=\Phg M_{3,u}$ and $\rk{\Rstar}{M_{1,u}}=\rk{\Rstar}{M_{3,u}}$
for all $\ee{1}{u}{m_r}$.
\end{lem}

\begin{pf}
Since $\Phg$ is a square matrix with 
$\det(\Phg)\neq0$ by \Lem{\ref{proc40}}, our assertion
follows from the preceding observation.
\end{pf}

\begin{lem}\lb{proc10}
With the notation above, suppose that $n^r_1 < n^{r-1}_{m_{r-1}}$. 
Put $p:=\max\{\ l\ |\ n^r_1+\kg_l\leq0,\ \ee{1}{l}{m_r}\ \}$ and
$\pcheck:=\max\{\ l\ |\ n^r_l+\kg_p\leq0,\ \ee{1}{l}{m_r}\ \}$, and
let $N_p:=0^p\opdot\lt(\dss{i=p+1}{m_r}\Rstar(\kg_i)\rt)\sset L$ be 
the $\Rstar$ module as in Section \ref{kiwadoi}.
Suppose further that
$\Rstar[\Xsub{1}\ddd\Xsub{\rg}]
\gradn{\Ima{\Phg}}{\leq -\kg_p}\sset N_p$. 
Then the following assertions hold,
where $(\ ,\ )$ denotes the canonical pairing of $L$ and $L^\vee$.
\begin{parenuma}
\item\lb{cond11}
$M_{1,\pcheck}\sset N_p$.
\item\lb{cond12}
$(v,\chg)=0$ for all $v\in M_{1,\pcheck}$ and $\chg\in M_{2,p}$.
\item\lb{cond31}
Given integers $u$ and $u'$ with $\ee{1}{u}{p}$, $\ee{1}{u'}{\pcheck}$,
there is an injective homomorphism
$\hg_{u,u'}:M_{2,u}\lra \Hom{\Rstar}{L/M_{1,u'}}{\Rstar}$
over $\Rstar[\tpose{\Xsub{1}}\ddd\tpose{\Xsub{\rg}}]$ 
under which the image of $\chg\in M_{2,u}$ is the homomorphism
$\hg_{u,u'}(\chg)$ defined by $\hg_{u,u'}(\chg)(\vbar)=(v,\chg)$ for 
$\vbar=v\ (\mod\ M_{1,u'})\in L/M_{1,u'}$.
\item\lb{cond13}
$\rk{\Rstar}{M_{1,\pcheck}}+\rk{\Rstar}{M_{2,p}}\leq\rk{\Rstar}{L}=m_r$.
\end{parenuma}
\end{lem}

\begin{pf}
Suppose $n^r_1 < r^{r-1}_{m_{r-1}}$. Then 
$n^r_1+\kg_1=n^r_1-n^{r-1}_{m_{r-1}}+1\leq0$ by
\Lem{\ref{proc44}}, so that $p$ and $\pcheck$ exist.
Let $\phgsub{i}$ be the $i$-th column of $\Phg$
for $\ee{1}{i}{m_r}$.
Since $\phgsub{1}$ is 
different from zero and its components lie in $x_r\Rstar$,
we see $p<m_r$, $\kg_p < \kg_{p+1}$.
Suppose further that $\Rstar[\Xsub{1}\ddd\Xsub{\rg}]
\gradn{\Ima{\Phg}}{\leq -\kg_p}\sset N_p$.
\begin{parenuma}
\item
For every $\ee{1}{l}{\pcheck}$, one has
$\deg(\Phg e^F_l)=n^r_l\leq -\kg_p$, in other words, 
$\Phg e^F_l\in\gradn{\Ima{\Phg}}{\leq -\kg_p}$.
Therefore 
$\Rstar[\Xsub{1}\ddd\Xsub{\rg}]\Phg e^F_l\in N_p$
for every $\ee{1}{l}{\pcheck}$ by our assumption. Hence 
$M_{1,\pcheck}\sset N_p$.
\item
Given $v\in M_{1,\pcheck}$ and $\chg\in M_{2,p}$,
$(\Xsub{i}v,\chg)=(v,\tpose{\Xsub{i}}\chg)$ for
all $\ee{1}{i}{\rg}$. It is therefore enough to verify that
$(v,e^{L\ast}_l)=0$ for all $v\in  M_{1,\pcheck}$ and $\ee{1}{l}{p}$.
But this holds since $M_{1,\pcheck}\sset N_p$ by \eqref{cond11}.
\item
The assertion can be verified with the use of \eqref{cond12}.
\item
By \eqref{cond31}, we can consider the homomorphism 
$\hg_{p,\pcheck}:M_{2,p}\lra \Hom{\Rstar}{L/M_{1,\pcheck}}{\Rstar}$.
Since $\hg_{p,\pcheck}$ is injective, 
$\rk{\Rstar}{M_{2,p}}\leq\rk{\Rstar}{L/M_{1,\pcheck}}
=\rk{\Rstar}{L}-\rk{\Rstar}{M_{1,\pcheck}}$, which proves our assertion.
\qed
\end{parenuma}
\killqed
\end{pf}

\par
Next we consider the modules, the matrices and the homomorphisms
which have appeared so far modulo $x_r$. 
Put $\Rbarstar:=\Rstar/x_r\Rstar=k$, 
$\Lbar:= L/x_rL=\dss{i=1}{m_r}k(\kg_i)$, 
$\Fbar:= F/x_rF=\dss{i=1}{m_r}k(-n^r_i)$, and
$\Lbar{}^\ast:=L^\vee/x_rL^\vee\cong\Hom{k}{L}{k}$. Denoting
$\Xsub{i}\ (\mod\ x_r)$ and $\Ysub{i}\ (\mod\ x_r)$ by 
$\Xbarsub{i}$ and $\Ybarsub{i}$ respectively, let further
\begin{align*}
&\Mbar_{2,u}=\smm{l=1}{u}
k[\tpose{\Xbarsub{1}}\ddd\tpose{\Xbarsub{\rg}}]
\ebar^{L\ast}_l\sset\Lbar{}^\ast,\\
&\Mbar_{3,u}=\smm{l=1}{u}
k[\Ybarsub{1}\ddd\Ybarsub{\rg}]\ebar^F_l
\sset\Fbar
\end{align*}
for $\ee{1}{u}{m_r}$. We have a commutative diagram
\begin{align*}
&\begin{CD}
\Lbar(-1) @>{\ \pgbar\ }>\sim>\Abar(-1)\\
 @VV \Xbarsub{i} V @VV \tm x_i V\\
\Lbar @>{\ \pgbar\ }>\sim>\Abar
\end{CD}\qquad.
\end{align*}
Besides, 
\begin{align}
&\Lbar\cong\Abar,\quad \Fbar\cong\Tor{1}{\Rstar}{A}{\Rbarstar}
\cong (0:_Ax_r)(-1),
\lb{eq48}\\
&\Lbar{}^\ast\cong\Hom{k}{\Abar}{k}
\cong\Ext{r-1}{\Rbar}{\Abar}{\Rbar}(-r+1),
\lb{eq17}\\
&\rk{k}{\Ext{r-1}{\Rbar}{\Abar}{\Rbar}}=m_r.
\lb{eq15}
\end{align}
Furthermore, if $A$ is \Gor,
\begin{equation}
\Mbar_{2,u}\cong\Mbar_{3,u}(n^r_{m_r})\quad\fallt\ \ee{1}{u}{m_r}
\lb{eq16}
\end{equation}
by \eqref{cond36} of \Rem{\ref{proc7}}.

\begin{lem}\lb{proc11}
With the same assumption and notation as in \Lem{\ref{proc10}},
we have $\rk{k}{\Mbar_{3,\pcheck}}+\rk{k}{\Mbar_{2,p}}\leq m_r$.
\end{lem}

\begin{pf}
Let $\seq{w}{1}{\dg}$ be elements of $M_{3,\pcheck}$ such that
$\seq{\wbar}{1}{\dg}\in\Mbar_{3,\pcheck}$ are linearly independent 
over $k$ with $\wbar_i:=w_i\ (\mod\ x_r)$. Then $\seq{w}{1}{\dg}$
are linearly independent over $\Rstar$, too. Since every element of
$\Mbar_{3,\pcheck}$ comes from that of $M_{3,\pcheck}$ modulo $x_r$,
we see $\rk{k}{\Mbar_{3,\pcheck}}\leq\rk{\Rstar}{M_{3,\pcheck}}$.
The same holds for $M_{2,p}$ and $\Mbar_{2,p}$.
Hence 
$\rk{k}{\Mbar_{3,\pcheck}}+\rk{k}{\Mbar_{2,p}}
\leq\rk{\Rstar}{M_{3,\pcheck}}+\rk{\Rstar}{M_{2,p}}$. 
On the other hand, $\rk{\Rstar}{M_{1,\pcheck}}=\rk{\Rstar}{M_{3,\pcheck}}$
by \Lem{\ref{proc12}}. Thus 
$\rk{k}{\Mbar_{3,\pcheck}}+\rk{k}{\Mbar_{2,p}}\leq m_r$ by 
\eqref{cond13} of \Lem{\ref{proc10}}.
\end{pf}

\begin{rem}\lb{proc21}
Observe that $\pg(\Phg e^F_l/x_r)\in 0:_A x_r$ by the argument
of the proof of \Prop{\ref{proc55}}.
Let $p$ and $\pcheck$ be as in \Lem{\ref{proc10}}.
Since $R$ and $A$ are commutative, we find by \Thm{\ref{proc27}} that 
$R\pg(\Phg e^F_l/x_r)\sset\pg(N_p)$
for all $\ee{1}{l}{\pcheck}$ if $\ch{k}=0$ and $\Gg$ is sufficiently
general. 
\end{rem}

\begin{lem}\lb{proc24}
Let $\seq{g}{1}{r}\in R$ be homogeneous polynomials of degrees 
$\seq{d}{1}{r}$ respectively which form an $R$-regular sequence
and let $I:=(\seq{g}{1}{r})\sset R$.
Suppose that $\Gg$ is chosen sufficiently generally.
\begin{parenuma}
\item\lb{cond34}
If $\ch{k}=0$ and $2d_i<d_1+\cdots+d_r-r+2$ for all $\ee{1}{i}{r}$,
then $\Ibar:=I+(x_r)/(x_r)\sset\Rbar$ is minimally generated by
$\seq{\gbar}{1}{r}$ over $R$, where $\gbar_i:=g_i\ (\mod\ (x_r))$.
\item\lb{cond35}
If $2d_i\geq d_1+\cdots+d_r-r+2$ for some $\ee{1}{i}{r}$, then
$R/I$ has the WLP.
\end{parenuma}
\end{lem}

\begin{pf}
For the case $r=3$, see \cite{JW}, too. Now, let us prove \eqref{cond34} first.
Suppose to the contrary. Rearranging $\seq{g}{1}{r}$ if necessary,
we may assume without loss of generality
that $\gbar_r\in(\seq{\gbar}{1}{r-1})$ for all sufficiently general $\Gg$.
Put $A':=R/(\seq{g}{1}{r-1})$ and $\Abar{}':=A'/x_rA'$. Then $\Abar=\Abar{}'$. 
Since $\seq{\gbar}{1}{r-1}$ form an $\Rbar$-regular sequence,
we see $0:_Ax_r\cong\Ext{r-1}{\Rbar}{\Abar}{\Rbar}(-r+1-\socd{A})
\cong \Abar{}'(-d_r+1)$ by \eqref{eq50}. The module $0:_Ax_r$ is therefore
generated over $R$ by $\gradn{0:_Ax_r}{\leq d_r-1}$. On the other hand,
$\socd{\Abar}=\socd{\Abar{}'}=d_1+\cdots+d_{r-1}-(r-1)$.
Hence we must have $d_r > d_1+\cdots+d_{r-1}-r+1$ by \eqref{cond38} of
\Prop{\ref{proc55}}, which contradicts our hypothesis.

\par
To prove \eqref{cond35}, we may assume with no loss of generality
that $2d_r\geq d_1+\cdots+d_r-r+2$, rearranging $\seq{g}{1}{r}$ if necessary.
Let the notation be the same as above. 
Since $\seq{g}{1}{r-1},x_r$ form an $R$-regular sequence, 
the $R$-module $A'$ is a free module over
$\Rstar$, so that there is a nondecreasing sequence $\seq{\kg}{1}{Q}$ such that
$A'\cong\dss{i=1}{Q}\Rstar(\kg_i)$. 
Here, $-\kg_1=\socd{\Abar{}'}=d_1+\cdots+d_{r-1}-(r-1)$.
We obtain a free resolution
\begin{equation}
0\lra\dss{i=1}{Q}\Rstar(-\kg'_i)\xra{\ \Phg\ }
\dss{i=1}{Q}\Rstar(\kg_i)\cong A'\lra A\lra 0
\lb{eq53}
\end{equation}
of $A:=R/I$ over $\Rstar$ with another nondecreasing sequence
$\kg'_1\ddd\kg'_{Q'}$. Let $\ghat_r$ denote $g_r\ (\mod\ (\seq{g}{1}{r-1}))$.
Since $\Ima{\Phg}=R\ghat_r$, we have $\kg'_1=d_r$.
Hence $\kg'_1+\kg_1=d_r+\kg_1\geq1$. 
This implies that the resolution \eqref{eq53} is minimal and that
$A$ has the WLP by the argument of the proof of
\Lem{\ref{proc5}}.
\end{pf}

\section{\mathversion{bold}The case where $r=3$ and
$0:_Ax_3$ is minimally generated by two elements}
\lb{Ibar}

\par
In this section, we consider homogeneous ideals $I$ in $R$ with $r=3$
defining \Gor\ graded rings $R/I$ of dimension zero such that
$0:_Ax_3$ is minimally generated by two elements.
For our purpose, it is necessary to have information on 
free resolutions of $\Ibar$ over $\Rbar$ in relation with free resolutions
of $I$ over $R$. Let us begin with some general properties.

\begin{lem}\lb{proc43}
Let $E$ be a \fgg\ $R$-module of $\depm{E}\geq r-s+1$ $(\ee{1}{s}{r})$,
$\Rbar:=R/(\seq{x}{s}{r})$, $\mfrakbar:=\mfrak/(\seq{x}{s}{r})$
and $E_{\Rbar}:=E/(\seq{x}{s}{r})E$.
Let further  
$\lgbar{}'_\bt:\Gbar{}'_\bt\lra E_{\Rbar}/\Hsupb{0}{\mfrakbar}{E_{\Rbar}}$ 
and $\lgbar{}''_\bt:\Gbar{}''_\bt\lra\Hsupb{0}{\mfrakbar}{E_{\Rbar}}$
be free resolutions of $E_{\Rbar}/\Hsupb{0}{\mfrakbar}{E_{\Rbar}}$ and 
$\Hsupb{0}{\mfrakbar}{E_{\Rbar}}$ 
respectively over $\Rbar$.
Suppose that the linear forms $\seq{x}{s}{r}$ form an $E$-regular
sequence.  
Then, there is a free resolution
$\lgbt:\Gbt\lra E$ of $E$ having the following properties.
\begin{parenumr}
\item\lb{cond14}
$G_j=G'_j\op G''_j$ for all $j\geq0$.
\item\lb{cond15}
Let $\vpg'_j:\Gsub{j}\lra\Gpsub{j}$ and $\vpg''_j:\Gsub{j}\lra\Gppsub{j}$
be the natural projections and put $\lg{}^{(1)}_j=\vpg'\cc\lg_j|_{\Gpsub{j}}$,
$\lg{}^{(3)}_j=\vpg''\cc\lg_j|_{\Gpsub{j}}$, 
$\lg{}^{(2)}_j=\vpg'\cc\lg_j|_{\Gppsub{j}}$,
and $\lg{}^{(4)}_j=\vpg''\cc\lg_j|_{\Gppsub{j}}$. 
Then $\lg{}^{(2)}_j\equiv0\ (\mod\ (\seq{x}{s}{r}))$ for all 
$j\geq1$.
\item\lb{cond16}
$(\Gbar{}'_\bt,\lgbar{}'_\bt)
=(G'_\bt\ot_R\Rbar,\lg{}^{(1)}_\bt\ot_R\Rbar)$ and
$(\Gbar{}''_\bt,\lgbar{}''_\bt)
=(G''_\bt\ot_R\Rbar,\lg{}^{(4)}_\bt\ot_R\Rbar)$.
\item\lb{cond17}
The map $\lg^{(3)}_1\ot_R\Rbar$ represents the element of 
$\Ext{1}{\Rbar}{E_{\Rbar}/\Hsupb{0}{\mfrakbar}{E_{\Rbar}}}
{\Hsupb{0}{\mfrakbar}{E_{\Rbar}}}$ 
corresponding to
the extension $E_{\Rbar}$ of $\Hsupb{0}{\mfrakbar}{E_{\Rbar}}$ by 
$E_{\Rbar}/\Hsupb{0}{\mfrakbar}{E_{\Rbar}}$.
\end{parenumr}
\end{lem}

\begin{pf}
Since $E_{\Rbar}$ lies in the middle of the exact sequence
\begin{equation}
0\lra\Hsupb{0}{\mfrakbar}{E_{\Rbar}}
\lra E_{\Rbar}\lra E_{\Rbar}/\Hsupb{0}{\mfrakbar}{E_{\Rbar}}\lra0,
\lb{eq30}
\end{equation}
there is a free resolution $\lgbarbt:\Gbarbt\lra E_{\Rbar}$ of $E_{\Rbar}$
such that 
\begin{equation*}
\Gbar_j=\Gbar{}'_j\op\Gbar{}''_j,\quad
\lgbar_j=\mat{\lgbar{}'_j}{0}{\lgbar{}'''_j}{\lgbar{}''_j}
\quad\fallt\quad j\geq0.
\end{equation*}
Let $G'_j$ (resp. $G''_j$) be the free $R$-module such that
$G'_j\ot_R\Rbar=\Gbar{}'_j$ (resp. $G''_j\ot_R\Rbar=\Gbar{}''_j$) and
$G_j:=G'_j\op G''_j$. Since the sequence $\seq{x}{s}{r}$ is 
$E$-regular by hypothesis, one can construct a free resolution 
$\lgbt:\Gbt\lra E$ satisfying all the 
conditions \eqref{cond15} -- \eqref{cond17}.
\end{pf}

\par
Given a module $E$ over $R$, let $\shmin{E}$ denote the number of
the minimal generators of $E$ over $R$. 

\begin{cor}\lb{proc45}
Let $I$ be a homogeneous ideal in $R$ 
such that the graded ring $A:=R/I$ is of dimension zero and let 
$\Abar:=A/x_rA$ and $\Ibar:=I+(x_r)/(x_r)$. Then
$\shmin{I}\leq\shmin{\Ibar}+\shmin{0:_A x_r}$. If further
$A$ is \Gor, then 
$\shmin{I}\leq\shmin{\Ibar}+\shmin{\Ext{r-1}{\Rbar}{\Abar}{\Rbar}}$.
\end{cor}

\begin{pf}
By \eqref{eq25}, $0:_A x_r\cong\Tor{1}{\Rstar}{A}{\Rbarstar}(1)
\cong\Hn{0}{\Itild}(1)$
and $\Itild/\Hn{0}{\Itild}\cong\Ibar$, where 
$\nfrak:=\mfrak/(x_r)\sset\Rbar$ and $\Itild:=I/x_rI$.
Our assertion therefore follows from \eqref{eq30} and \eqref{eq50} since 
$\shmin{I}=\shmin{\Itild}$.
\end{pf}

\par
In the rest of this section, unless otherwise specified,
we assume that $r=3$ and that $I$ is a
homogeneous ideal in $R$ defining a \Gor\ graded ring 
$A:=R/I$ of dimension zero with the property $\shmin{0:_Ax_3}=2$.
Notice that the condition $\shmin{0:_Ax_3}=2$ is equivalent to
$\shmin{\Ibar}=3$ by \eqref{eq50} in our case, where
$\Ibar:=I+(x_3)/(x_3)\sset\Rbar:=R/(x_3)=k[x_1,x_2]$.

\begin{ex}\lb{proc63}
There is indeed a \Gor\ graded ring with the property mentioned above which
is not a complete intersection. 
Put $g_1:=-y_1^5$, $g_2:=y_1^2y_2^3$, $g_3:=y_1^3y_3^3-y_2^6$, $g_4:=y_2^3y_3^4$,
$g_5:=-y_3^7$, and $I:=(g_1,g_2,g_3,g_4,g_5)\sset R=k[y_1,y_2,y_3]$.
For each $\ee{1}{i}{5}$, 
the polynomial $g_i$ is equal to the Pfaffian of $\Lg\bnmt{i}{i}$ up to sign, 
where
\begin{equation*}
\Lg=
\begin{bmatrix}
0&0&0&y_3^3&y_2^3\\
0&0&y_3^4&y_2^3&y_1^3\\
0&-y_3^4&0&y_1^2&0\\
-y_3^3&-y_2^3&-y_1^2&0&0\\
-y_2^3&-y_1^3&0&0&0
\end{bmatrix}
\end{equation*}
and $\Lg\bnmt{i}{i}$ denotes
the alternating matrix obtained from $\Lg$ by deleting its $i$-th row
and column. Hence the ring $R/I$ is \Gor\ by 
\cite[\Thm{2.1}]{BE}. Let $a$ and $b$ be sufficiently general 
elements of $k$, and put $x_3:=y_3+ay_1+by_2$, $\gbar_i:=g_i(y_1,y_2,-ay_1-by_2)$
$(\ee{1}{i}{5})$ and $\Ibar:=(\gbar_1,\gbar_2,\gbar_3,\gbar_4,\gbar_5)
\sset\Rbar=R/(x_3)$. By direct computations one can verify that
$\Ibar=(\gbar_1,\gbar_2,\gbar_3)$. Since $J:=(g_1,g_2,g_3)\sset R$ is generated by the
maximal minors of 
$\begin{bmatrix}
y_3^3&y_2^3&y_1^2\\
y_2^3&y_1^3&0
\end{bmatrix}$, one can compute $\dimk{\gradn{J}{l}}$ and
$\dimk{\gradn{\Ibar}{l}}$ for all $l\in\Zbf$ 
and find that $\BR{J}=\BRp{\Ibar}=(5;5,6,6,7,7)$. Since $\deg(g_4)=\deg(g_5)=7$,
it follows that $n^3_1=7=n^2_5$. Hence $R/I$ has the WLP by \Lem{\ref{proc5}}.
In fact $\BR{I}=(5;5,6,6,7,7;7,7,8,8,8,8,9,9,9,9,9,10,10,10,10,11,11,11,12,12,13)$.
\end{ex}

\begin{lem}\lb{proc35}
Suppose that $\ch{k}=0$, $\shmin{0:_Ax_3}=2$, and that
$\Gg$ is chosen sufficiently generally.
Let $(a;n^2_1\ddd n^2_a;n^3_1\ddd n^3_b)$ be the basic sequence of $I$. 
Then, there are integers 
$n,d_1,d_2,d_3$ with $0<d_1\leq d_2\leq d_3$, $d_2<n\leq n^3_1+1$,
$n\leq n^2_a+1$ such that 
$\Ibar$ has a minimal free resolution of the form
\begin{equation}
0\lra\Rbar(-n)\op\Rbar(-n^2_a-1)\xra{\ \Psg\ }
\Rbar(-d_1)\op\Rbar(-d_2)\op\Rbar(-d_3)
\lra\Ibar\lra0
\lb{eq39}
\end{equation}
over $\Rbar$.
\end{lem}

\begin{pf}
The basic sequence of $\Ibar$ is $(a;n^2_1\ddd n^2_a)$, so that 
$\Ibar$ has a free resolution of the form
\begin{equation}
0\lra\dss{i=1}{a}\Rbar(-n^2_i-1)\lra
\Rbar(-a)\op\dss{i=1}{a}\Rbar(-n^2_i)
\lra\Ibar\lra0
\lb{eq54}
\end{equation}
(see e.g. \cite[pp. 10 -- 13]{A8}).
Since $\shmin{\Ibar}=3$ by hypotheses, after elimination, 
the above resolution reduces to the form \eqref{eq39},
where $d_1=a$, $d_2=n^2_1$, $d_3=n^2_i$ and $n=n^2_j+1$ for some  $i,j$ 
with $\ne{1}{i}{a}$, $\en{i-1}{j}{a}$.
It remains to show that $n\leq n^3_1+1$.
Since $\socd{A}=n^3_b-1$ and 
$0:_Ax_3\cong\Ext{2}{\Rbar}{\Abar}{\Rbar}(-2-\socd{A})$, it follow from 
\eqref{eq39} that the module
$0:_Ax_3$ is generated by $\gradn{0:_Ax_3}{\leq n^3_b+1-n}$ over $R$.
By \eqref{cond26} of \Prop{\ref{proc55}} and \Lem{\ref{proc37}}, we must have 
$(n^3_b+1-n+1)+(-n^3_b+n^3_1) > 0$. Hence $n^3_1+1\geq n$.
\end{pf}

\par
For a sequence $B=(\ag;\ngbar^2;\ngbar^3)
=(\ag;\ng^2_1\ddd\ng^2_\ag;\ng^3_1\ddd\ng^3_\bg)$
of integers, put
\begin{equation*}
h_B(\mg):=\bnmm{\mg-\ag+2}{2}+\smm{l=1}{\ag}\bnmm{\mg-\ng^2_l+1}{1}+
\smm{l=1}{\bg}\bnmm{\mg-\ng^3_l}{0}
\end{equation*}
for $\mg\in\Zbf$, where
\begin{equation*}
\bnmm{\mg}{\ng}:=\frac{\mg!}{(\mg-\ng)!\ng!}\quad
\text{if $\mg\geq\ng\geq0$ \ and}\quad
\bnmm{\mg}{\ng}:=0\quad
\text{otherwise.}
\end{equation*}
Given sequences $B$ as above and
$B'=(\ag';\ngbar{}'^2;\ngbar{}'^3)
=(\ag';\ng'{}^2_1\ddd\ng'{}^2_{\ag'};\ng'{}^3_1\ddd\ng'{}^3_{\bg'})$,
we will write $B\sim B'$ to mean that $h_B(\mg)=h_{B'}(\mg)$
for all $\mg\in\Zbf$.
We further define a relation $B \eqp B'$ to mean that 
$\ag=\ag'$ and that $\ngbar^i = \ngbar'{}^i$ up to permutation for 
$i=2,3$.
Given a sequence $(l_{j1}\ddd l_{ji_j})$ $(\ee{1}{j}{m})$, we will
denote by $((l_{j1}\ddd l_{ji_j})_{\ee{1}{j}{m}})$ the sequence
$(l_{11}\ddd l_{1i_1},l_{21}\ddd l_{2i_2},\cdots, l_{m1}\ddd l_{mi_m})$.

\begin{lem}\lb{proc32}
Let $g_1,g_2$ and $g_3$ be homogeneous polynomials in 
$R$ of degrees $d_1,d_2$ and $d_3$ $(1\leq d_1\leq d_2\leq d_3)$
respectively which form an $R$-regular sequence.
Let further $I:=(g_1,g_2,g_3)\sset R$ be the homogeneous 
ideal generated by them. Then
$\BR{I}\sim B^0$, where
\begin{align*}
\ngbar^{02} &:= (d_2,d_2+1\ddd d_2+d_1-1),\\
\ngbar^{03} &:= ((d_3+j,d_3+1+j\ddd d_3+d_1-1+j)_{\ee{0}{j}{d_2-1}})
\end{align*}
and $B^0:=(d_1;\ngbar^{02};\ngbar^{03})$.
\end{lem}

\begin{pf}
It follows from the minimal free resolution of $I$ over $R$
that
\begin{align}
h_{\BR{I}}(\mg)&=\dimk{[I]_{\mg}}
\notag\\
&\hspace{-1em}=\bnmm{\mg-d_1+2}{2}+\bnmm{\mg-d_2+2}{2}+\bnmm{\mg-d_3+2}{2}
\notag\\
&\hspace{1em}-\bnmm{\mg-d_1-d_2+2}{2}-\bnmm{\mg-d_2-d_3+2}{2}
-\bnmm{\mg-d_3-d_1+2}{2}
\notag\\
&\hspace{1em}+\bnmm{\mg-d_1-d_2-d_3+2}{2}
\lb{eq40}\\
&\hspace{-1em}=\bnmm{\mg-d_1+2}{2}+\smm{i=0}{d_1-1}\bnmm{\mg-d_2-i+1}{1}
+\smm{i=0}{d_1-1}\smm{j=0}{d_2-1}\bnmm{\mg-d_3-i-j}{0}
\notag\\
&=h_{B^0}(\mg)
\notag
\end{align}
for all $\mg\in\Zbf$. This implies our assertion.
\end{pf}

For the sake of completeness, we will include a proof for the height 
three Artinian complete intersection
case in our argument of the main \Thm{\ref{proc61}}, so that we need
the following lemma which may follow from the results of \cite{JW} and 
\cite{HMNJW}.

\begin{lem}\lb{proc57}
With the notation of \Lem{\ref{proc32}}, let
$(a;n^2_1\ddd n^2_a;n^3_1\ddd n^3_b)$ be the basic sequence of $I$ and 
$A:=R/I$. 
Let further $n$ be the integer stated in \Lem{\ref{proc35}}.
Suppose that $\Gg$ is chosen sufficiently generally and
that $\ch{k}=0$.
\begin{parenuma}
\item\lb{cond27}
If $d_3\geq d_1+d_2-1$, then $\BR{I}\eqp B^0$. In particular,
$n^2_a=d_2+d_1-1\leq d_3=n^3_1$.
\item\lb{cond28}
If $d_3 < d_1+d_2-1$ and $n^3_1 < n^2_a$, then $n^2_a+1-n\geq 2$.
\end{parenuma}
\end{lem}

\begin{pf}
Suppose first that $d_3\geq d_1+d_2-1$. Then $A$ has the WLP by \eqref{cond34}
of \Lem{\ref{proc24}}, so that $n^3_1\geq n^2_a$ by \Lem{\ref{proc5}}.
Since the smallest term of $\ngbar^{03}$ is not less than the largest
term of $\ngbar^{02}$, we find by the equality $h_{\BR{I}}(\mg)=h_{B^0}(\mg)$
$(\mg\in\Zbf)$ that $\BR{I}\eqp B^0$.

\par
Next assume that $d_3 < d_1+d_2-1$ and that $n^3_1 < n^2_a$.
By \Lem{\ref{proc35}}, we have $n\leq n^3_1+1\leq n^2_a<n^2_a+1$, so that
$n^2_a+1-n\geq 1$. Since \eqref{eq54} reduces to \eqref{eq39} by eliminations, 
the transpose of the matrix $\Psg$ appearing in \eqref{eq39} must be of the form 
$\tpose{\Psg}=
\begin{pmatrix}
\gbar_{11}&\gbar_{12}&\gbar_{13}\\
\gbar_{21}&\gbar_{22}&\gbar_{23}
\end{pmatrix}$ with (a) $\gbar_{13}=0$, $\deg(\gbar_{23})=1$ or 
(b) $\deg(\gbar_{23})\geq 2$. 
Let $\gbar_i:=g_i\ (\mod\ (x_3)) \in \Rbar$.
In the case (a), $\gbar_1=\gbar_{23}\gbar_{12}$, 
$\gbar_2=\gbar_{23}\gbar_{11}$ up to multiplication by a constant 
respectively. But $\gbar_1,\gbar_2$ must be relatively prime, 
since the polynomials $g_1,g_2,g_3$ form an $R$-regular sequence and
$\Gg$ is sufficiently general. Hence only the case (b) occurs.
By \eqref{eq50} and \eqref{eq48}, 
\begin{equation}
\dss{i=1}{b}k(-n^3_i+1)\cong 0:_Ax_3\cong\Hom{k}{\Abar}{k}(-n^3_b+1)
\cong\Ext{2}{\Rbar}{\Abar}{\Rbar}(-n^3_b-1),
\lb{eq41}
\end{equation}
and $n^3_1-1=n^3_b-n^2_a$ by \Cor{\ref{proc39}}. 
Now consider the case $n^2_a+1-n=1$. First, 
$\dimk{\gradn{0:_Ax_3}{n^3_1}}=3$ by \eqref{eq41} and the form of $\Psg$.
Secondly, $n=n^3_1+1=n^2_a$, so that
$2n^3_1=n^3_b=d_1+d_2+d_3-2$.
Observe that $n^3_1-2-d_i\geq -2$ $(i=1,2,3)$ and that 
$n^3_1-d_i-d_j\leq -2$ $(1\leq i<j\leq 3)$.
The map $\tm x_3:[A]_{n^3_1-2}\lra [A]_{n^3_1-1}$ is injective 
by the argument of the proof of \Lem{\ref{proc5}},
$\gradn{0:_Ax_3}{n^3_1}=[\Hom{k}{\Abar}{k}]_{-n^3_1+1}$ 
by \eqref{eq41}, and 
$\bnmmt{\mg+1}{2}-\bnmmt{\mg}{2}=\mg$ for all $\mg\geq0$. Hence,
\begin{align*}
\dimk{\gradn{0:_Ax_3}{n^3_1}}
&=\dimk{[\Abar]_{n^3_1-1}}=\dimk{\gradn{A}{n^3_1-1}}-\dimk{\gradn{A}{n^3_1-2}}\\
&\hspace{-2em}=\dimk{\gradn{R}{n^3_1-1}}-\dimk{\gradn{R}{n^3_1-2}}
-\dimk{\gradn{I}{n^3_1-1}}+\dimk{\gradn{I}{n^3_1-2}}=2
\end{align*}
by \eqref{eq40}, which is a contradiction. Thus $n^2_a+1-n\geq 2$.
\end{pf}

\begin{lem}\lb{proc58}
Notation being as in \Lem{\ref{proc35}}, let further
$\Phg$, $\pg$, $\thgsub{ijl}$ and $\Xsub{i}$ be as 
in Sections \ref{mat} and \ref{gor}.
\begin{parenuma}
\item\lb{cond29}
If $n^3_1 < n^2_a$, 
$I=(g_1,g_2,g_3)$ with homogeneous polynomials $g_1,g_2$ and $g_3$
of degrees $d_1,d_2$ and $d_3$ $(1\leq d_1\leq d_2\leq d_3)$ respectively
and $d_3 < d_1+d_2-1$, then 
$\gradn{\Ima{\Phg}}{n^3_1+1}=
\smm{i=1}{2}\Xsub{i}\gradn{\Ima{\Phg}}{n^3_1}
+x_3\gradn{\Ima{\Phg}}{n^3_1}$.
\item\lb{cond30}
If $\shmin{I}\neq3$, then $\shmin{I}=5$ and 
$\thgsub{12l}\in\smm{i=1}{2}\Xsub{i}\gradn{\Ima{\Phg}}{n^3_b-n^3_l+1}+
\linebreak
x_3\gradn{\Ima{\Phg}}{n^3_b-n^3_l+1}$ for all
$\ee{1}{l}{b}$.
\end{parenuma}
\end{lem}

\begin{pf}
\no{1}
We have $n^2_a+1-n\geq 2$ by \eqref{cond28} of \Lem{\ref{proc57}} and
$n^3_b-n^2_a=n^3_1-1$ by \Cor{\ref{proc39}}, so that 
$\gradn{0:_Ax_3}{n^3_1}=\smm{i=1}{2}x_i\gradn{0:_Ax_3}{n^3_1-1}$ by
\eqref{eq41}. As in the proof of \Prop{\ref{proc55}}, we see that 
\begin{equation*}
\gradn{\Imasup{\Rstar}{\Phg/x_r}}{n^3_1}\sset
\smm{i=1}{2}\Xsub{i}\gradn{\Imasup{\Rstar}{\Phg/x_r}}{n^3_1-1}
+\gradn{\Ima{\Phg}}{n^3_1}, 
\end{equation*}
from which follows our assertion.
\indno{2}
Let $\{f^1_1,\ f^2_1\ddd f^2_a,\ f^3_1\ddd f^3_b\}$ be 
a pseudo \Wei\ basis of $I$ with respect to $x_1,x_2,x_3$ and let
$\lgbar{}'_\bt:\Gbar{}'_\bt\lra \Itild/\Hn{0}{\Itild}=\Ibar$ and 
$\lgbar{}''_\bt:\Gbar{}''_\bt\lra\Hn{0}{\Itild}$
be free resolutions of $\Ibar$ and 
$\Hn{0}{\Itild}$ respectively over $\Rbar$ such that
$\lgbar{}'_0=(\fbar{}^1_1\ \fbar{}^2_1\ \ldots\ \fbar{}^2_a)$,
$\lgbar{}''_0=(\ftild{}^3_1\ \ldots\ \ftild{}^3_b)$.
Here, we can choose $\lgbar{}'_\bt$ so that there is no unit in the first row
of $\lgbar{}'_1$. Then, there is a free resolution
$\lgbt:\Gbt\lra I$ of $I$ having the properties stated in
\Lem{\ref{proc43}}.
Since $\shmin{\Ibar}=3$ by assumption, the rank of 
$\lgbar{}'_1\ (\mod\ \mfrakbar)$ over $k$ is $a-2$.
After a suitable column transformation of $\lg_1|_{G'_1}$ and
a change of numbering of $f^2_1\ddd f^2_a$, 
if necessary, the relations given by the columns of $\lg_1|_{G'_1}$
are of the forms
\begin{align*}
&h_{1i}f^1_1+\smm{j=1}{2}h_{j+1\,i}f^2_j+
\smm{l=1}{b}h_{l+a+1\,i}f^3_l=0 \quad(\ee{1}{i}{2}),\\
&h_{1i}f^1_1+\smm{j=1}{2}h_{j+1\,i}f^2_j+
f^2_i+\smm{l=1}{b}h_{l+a+1\,i}f^3_l=0 \quad(\ee{3}{i}{a}).
\end{align*}
Put $f'{}^2_i:=f^2_i+\smm{l=1}{b}h_{l+a+1\,i}f^3_l$ for $\ee{3}{i}{a}$,
$f_1:=f^1_1$, $f_2:=f^2_1$ and $f_3:=f^2_2$.
Then $P:=\{f^1_1,\ f^2_1,f^2_2,f'{}^2_3\ddd f'{}^2_a,\ f^3_1\ddd f^3_b\}$
is a pseudo \Wei\ basis of $I$ with respect to $x_1,x_2,x_3$ by 
\Lem{\ref{proc51}}. We may further assume by \Lem{\ref{proc53}} that
$f^3_l=\smm{i=1}{b}\phgsub{il}\gcheck_i$ for all
$\ee{1}{l}{b}$ with the notation there. 
Now the sets $P$, $P_1:=P-\{f^3_1\ddd f^3_b\}$ and the polynomials 
$f_1,f_2,f_3$ satisfy the
conditions \eqref{cond23} and \eqref{cond24} required in \Lem{\ref{proc54}}.
To verify \eqref{eq38}, recall the exact sequence
\begin{equation*}
0\lra\Hn{0}{\Itild}\lra\Itild\lra\Ibar\lra0.
\end{equation*}
Since $\shmin{\Ibar}=3$ by hypotheses, we see 
$3\leq\shmin{\Itild}=\shmin{I}\leq 3+2=5$ by
\Cor{\ref{proc45}}. On the other hand, $\shmin{I}$ must be odd since $A$
is \Gor. Besides, $\shmin{I}\neq3$ again by hypotheses.
Hence $\shmin{\Itild}=5$. This implies that the sequence
\begin{equation*}
0\lra\Hn{0}{\Itild}/\nfrak\Hn{0}{\Itild}\lra
\Itild/\nfrak\Itild\lra\Ibar/\nfrak\Ibar\lra0
\end{equation*}
is also exact. Hence the condition \eqref{eq38} holds and our assertion
follows from \Lem{\ref{proc54}}.
\end{pf}

From now on through the end of this section,
let the  notation and the assumption be as in \Lem{\ref{proc35}}.
Let further $\gbar_1,\ \gbar_2$ and $\gbar_3$ be homogeneous elements of $\Ibar$
with $g_i\in I$ and $\deg(\gbar_i)=d_i$ $(i=1,2,3)$ giving a system of 
minimal generators of $\Ibar$. 
Taking the dual of the minimal free resolution described in \Lem{\ref{proc35}},
we get an exact sequence
\begin{align}
0\lra&\Rbar(n^2_a+1)/\Imasup{\Rbar}{\tpose{\Psg}}\cap\Rbar(n^2_a+1)
\xra{\ \ig\ }\Ext{2}{\Rbar}{\Abar}{\Rbar}
\notag\\
\lra&\Rbar(n)/(\gbar_{11},\gbar_{12},\gbar_{13})\lra0,
\lb{eq22}
\end{align}
where $\tpose{\Psg}=
\begin{pmatrix}
\gbar_{11}&\gbar_{12}&\gbar_{13}\\
\gbar_{21}&\gbar_{22}&\gbar_{23}
\end{pmatrix}$.

\begin{lem}\lb{proc36}
Suppose that $n^3_1<n^2_a$.
If $\gbar_1$ and $\gbar_2$ form an $\Rbar$-regular sequence or 
$\gbar_{13}\neq0$, then
\begin{equation*}
\rk{k}{\Rbar(n^2_a+1)/\Imasup{\Rbar}{\tpose{\Psg}}\cap\Rbar(n^2_a+1)}
>\rk{k}{\Rbar(n)/(\gbar_{11},\gbar_{12},\gbar_{13})}.
\end{equation*}
\end{lem}

\begin{pf}
Let 
\begin{align*}
0&\lra\Rbar(-\sg_1)\op\Rbar(-\sg_2)\xra{\ \Psg'\ }
\Rbar(d_1-n)\op\Rbar(d_2-n)\op\Rbar(d_3-n)\\
&\xra{\ (\gbar_{11},\gbar_{12},\gbar_{13})\ }\Rbar\lra
\Rbar/(\gbar_{11},\gbar_{12},\gbar_{13})\lra0
\end{align*}
be a free resolution of $\Rbar/(\gbar_{11},\gbar_{12},\gbar_{13})$, where
$\sg_1\leq\sg_2$ and 
$\tpose{\Psg'}=
\begin{pmatrix}
\gbar'_{11}&\gbar'_{12}&\gbar'_{13}\\
\gbar'_{21}&\gbar'_{22}&\gbar'_{23}
\end{pmatrix}$.

\par
If $\sg_2 < n-d_1$, then $\sg_1 < n-d_1$, $\deg(\gbar'_{11})<0$,
$\deg(\gbar'_{21})<0$, and $\gbar'_{11}=\gbar'_{21}=0$. 
But $\rank{\Psg'}=2$, so that $\gbar_{12}=\gbar_{13}=0$.
Since $\rank{\Psg}=2$, this implies that $\gbar_{11}\neq0$ and that
$\gbar_1=0$, which is impossible.
If $\sg_1 < n-d_2$, then $\gbar'_{11}=\gbar'_{12}=0$, and
$\gbar'_{13}\neq0$, so that $\gbar_{13}=0$ and $\gbar'_{13}$
must lie in $k^\ast$. Consequently,  
$\tpose{\Psg}=
\begin{pmatrix}
\gbar_{11}&\gbar_{12}&0\\
\gbar_{21}&\gbar_{22}&\gbar_{23}
\end{pmatrix}$, and
$\gbar_1=\gbar_{23}\gbar_{12}$, $\gbar_2=\gbar_{23}\gbar_{11}$ and
$\gbar_3=\gbar_{11}\gbar_{22}-\gbar_{12}\gbar_{21}$
up to multiplication by elements of $k^\ast$.
But that cannot happen when $\gbar_{13}\neq0$.
In the case where $(\gbar_1,\gbar_2)$ is a complete intersection, 
we see
$\gbar_{23}\in k^\ast$. Hence $\shmin{\Ibar}$ must be two,
contradicting our assumption.

\par
The above observation implies that $\sg_2\geq n-d_1$ and that
$\sg_1 \geq n-d_2$. Put $\tg_i:=\sg_i+n^2_a+1-n$ for
$i=1,2$. We obtain an exact sequence of the form
\begin{align*}
0&\lra\Rbar(-\tg_1-\tg_2)(n^2_a+1)\lra(\Rbar(-\tg_1)\op\Rbar(-\tg_2))(n^2_a+1)\\
&\xra{\ (\gbar_{21},\gbar_{22},\gbar_{23})\Psg'\ }
\Rbar(n^2_a+1)\lra
\Rbar(n^2_a+1)/\Imasup{\Rbar}{\tpose{\Psg}}\cap\Rbar(n^2_a+1)\lra0.
\end{align*}
Since $n^2_a>n^3_1$ by our hypothesis, we see $n^2_a+1-n>0$ by
\Lem{\ref{proc35}}. Hence
\begin{align*}
&\rk{k}{\Rbar(n^2_a+1)/\Imasup{\Rbar}{\tpose{\Psg}}\cap\Rbar(n^2_a+1)}\\
&\hskip3em=\tg_1\tg_2 > \sg_1\sg_2 \geq (n-d_2)(n-d_1)
\geq \rk{k}{\Rbar(n)/(\gbar_{11},\gbar_{12},\gbar_{13})},
\end{align*}
as desired.
\end{pf}

\begin{lem}\lb{proc59}
Put $d_4=n^3_1$ and $d_5=n^3_b-n+2$.
If $\deg(\gcd(g_1,g_2))>0$ and $\gbar_{13}=0$, 
there is a  minimal free resolution of $I$ over $R$ of the form
\begin{equation}
0\lra R(-n^3_b-2)\lra\dss{i=1}{5}R(-n^3_b-2+d_i)\xra{\ \Lg\ }
\dss{i=1}{5}R(-d_i)\xra{\ \lg_0\ } I\lra0,
\lb{eq42}
\end{equation}
where
\begin{equation}
\Lg=
\begin{bmatrix}
0&\ast&\ast&h_3&h_1\\
\ast&0&\ast&h_4&h_2\\
\ast&\ast&0&h_5&0\\
-h_3&-h_4&-h_5&0&0\\
-h_1&-h_2&0&0&0
\end{bmatrix},
\lb{eq43}
\end{equation}
$g_1=-h_2h_5$, $g_2=h_1h_5$, $g_3=h_2h_3-h_1h_4$,
$h_5=\gcd(g_1,g_2)$, and
$d_1\leq d_2\leq d_3 \leq d_4 \leq d_5$.
\end{lem}

\begin{pf}
Suppose that $\deg(\gcd(g_1,g_2))>0$ and that $\gbar_{13}=0$.
Then $g_1=-h_2\gcd(g_1,g_2)$ and $g_2=h_1\gcd(g_1,g_2)$ with
relatively prime homogeneous polynomials $h_1,h_2\in R$ of positive
degree. Since $\Gg$ is sufficiently general,
$\hbar_1$ and $\hbar_2$ are also relatively prime. The syzygy module
\begin{equation*}
\Syz{\Rbar}{1}{\gbar_1,\gbar_2}:=
\{\ (\hbar{}',\hbar{}'')\in\Rbar(-d_1)\op\Rbar(-d_2)\ |
\ \hbar{}'\gbar_1+\hbar{}''\gbar_2=0\ \}
\end{equation*}
is therefore generated by
$(\hbar_1,\hbar_2)$. On the other hand, it is also generated by
$(\gbar_{11},\gbar_{12})$ since $\gbar_{13}=0$ by our hypothesis.
Hence $(\hbar_1,\hbar_2)$ coincides with $(\gbar_{11},\gbar_{12})$
up to multiplication by a constant and $\deg(\hbar_i)=\deg(\gbar_{1i})$
for $i=1,2$. In consequence, $\deg(g_1h_1)=\deg(g_2h_2)=n$.
Besides, $\gbar_{23}=\gcd(\gbar_1,\gbar_2)$ up to multiplication by a constant.

\par
By what we have seen in the proof of \Lems{\ref{proc40}} and by \eqref{eq50}, 
\begin{equation}
\dss{l=1}{b}k(-n^3_l)\cong\dss{l=1}{b}k\ftild^3_l=\Hsupb{0}{\nfrak}{\Itild}
\cong\Ext{2}{\Rbar}{\Rbar/\Ibar}{\Rbar}(-n^3_b-2),
\lb{eq44}
\end{equation}
since $\Tor{1}{\Rstar}{A}{\Rbarstar}\cong (0:_Ax_3)(-1)$.
By \eqref{eq39} and \eqref{eq44}, 
$\Hn{0}{\Itild}$ is generated
minimally by two elements. Let $g_4$ and $g_5$ be elements of $\supang{I}{3}$
with $\deg(g_4)\leq\deg(g_5)$ such that $\gtild_4$ and $\gtild_5$ generate
$\Hn{0}{\Itild}$ minimally over $\Rbar$. Here $\deg(g_4)=n^3_1$ 
and $\deg(g_5)=n^3_b-n+2$ by \eqref{eq44}.
Since $\deg(\gcd(g_1,g_2))>0$ by
our hypothesis, $\shmin{I}\neq3$. Hence
$\shmin{I}=5$ by \eqref{cond30} of \Lem{\ref{proc58}} 
and $I$ is generated minimally by $g_i$ $(\ee{1}{i}{5})$ over $R$.
Notice that $d_1\leq d_2\leq\deg(g_i)$ for all $i=3,4,5$ by \Lem{\ref{proc35}}.
With the use of \Lem{\ref{proc43}}, we find that there is a surjection
\begin{equation*}
R(-n)\op R(-n^2_a-1)\op \dss{j=1}{3}R(-n^3_1-\deg(\gbar_{2j}))
\thra\Syz{R}{1}{I}:=\Kersup{R}{\lg_0}.
\end{equation*}
Since $\deg(\gbar_{2j})>0$ for all $i=1,2,3$, we have
$n=\min\{\ l\ |\ \gradn{\Syz{R}{1}{I}}{l}\neq0\ \}$ by \Lem{\ref{proc35}}.

\par
On the other hand, since $A$ is \Gor,
there is a minimal free resolution of $I$ of the form
\begin{equation}
0\lra R(-\bg)\lra\dss{i=1}{5}R(-\bg+\ag_i)\xra{\ \Lg\ }
\dss{i=1}{5}R(-\ag_i)\lra I\lra0
\lb{eq45}
\end{equation}
with an alternating matrix $\Lg$, where 
$0<\ag_1=d_1\leq d_2=\ag_2\leq\cdots\leq\ag_5$
(see \cite[\Thm{2.1}]{BE}).
Since $\tpose{(h_1,h_2,0,0,0)}$ lies in $\gradn{\Syz{R}{1}{I}}{n}$,
we may assume with no loss of generality that $\Lg$ is of the
form \eqref{eq43}.
Put $g'_1:=-h_2h_5$, $g'_2:=h_1h_5$, $g'_3:=h_2h_3-h_1h_4$ and 
$J:=(g'_1,g'_2,g'_3)$. Then  
$\Pfaff(\Lg\bnmt{i}{i})$ equals $g'_i$ up to multiplication by
a constant respectively for $i=1,2,3$, where $\Lg\bnmt{i}{j}$ denotes
the square matrix obtained from $\Lg$ by deleting its $i$-th row
and $j$-th column. Hence $J\sset I$. 
Notice that $\deg(g'_i)=d_i$ for $i=1,2$ and that
$\deg(h_5)=\deg(g'_1)-\deg(h_2)=d_1-\deg(h_2)=\deg(\gcd(g_1,g_2))>0$.
If $\hei{J}\leq1$,
then there is an irreducible polynomial $h'$ dividing $h_5$ such that
$(h_1,h_2)$ and $(h_3,h_4)$ $(\mod\ (h'))$ are 
linearly dependent over $R/(h')$. In this case, $\hei{I}$ cannot be three 
since $I$ is generated by $\Pfaff(\Lg\bnmt{i}{i})$ $(\ee{1}{i}{5})$.
Hence $\hei{J}=2$. In consequence, $\hei{\Jbar}=2$ and $\Syz{\Rbar}{1}{\Jbar}$
is generated by $\tpose{(\hbar_1,\hbar_2,0)}$
and $\tpose{(\hbar_3,\hbar_4,\hbar_5)}$, so that 
$\Ibar\ni\gbar{}'_3\notin(\gbar{}'_1,\gbar{}'_2)$. 
If $\deg(g_3)=d_3>\ag_3=\deg(\gbar{}'_3)$, then 
$(\gbar{}'_1,\gbar{}'_2)$ must coincide with $(\gbar_1,\gbar_2)$, 
since $(\gbar{}'_1,\gbar{}'_2)\sset\Ibar=(\gbar_1,\gbar_2,\gbar_3)$
and $\deg(\gbar{}'_i)=d_i$ $(i=1,2)$.
In addition, $\gbar{}'_3$ must be contained in
$(\gbar_1,\gbar_2)=(\gbar{}'_1,\gbar{}'_2)$, which is a contradiction.
This implies that
$\ag_3\geq d_3$. But $d_3$ coincides with some 
$\ag_i$ with $i\geq3$. Hence $\ag_3=d_3$.
Moreover $\ag_4=\deg(g_4)=n^3_1$ and $\ag_5=\deg(g_5)=n^3_b-n+2$. 
Since $\Jbar\sset\Ibar$, $\deg(\gbar{}'_i)=d_i$ $(i=1,2,3)$, and
$\Jbar$ is generated minimally by $\gbar{}'_1,\gbar{}'_2,\gbar{}'_3$ over 
$\Rbar$, one sees that $\Jbar=\Ibar$. 
One may therefore assume from the first that
$\gbar{}'_i=\gbar_i$, $g'_i=g_i$ for $i=1,2,3$. 
Finally, $\socd{A}=n^3_b-1=\bg-3$.
Thus the free resolution \eqref{eq45} is of the desired form \eqref{eq42}.
\end{pf}

\begin{lem}\lb{proc60}
Suppose that $\deg(\gcd(g_1,g_2))>0$ and that $\gbar_{13}=0$.
With the notation of \Lem{\ref{proc59}}, put $J:=(g_1,g_2,g_3)\sset I$.
Then $R/J\cong\dss{i=1}{b}\Rstar(-n^3_b+n^3_i)$ and 
$I/J\cong\Ext{2}{R}{R/J}{R}(-\socd{A}-3)\cong\dss{i=1}{b}\Rstar(-n^3_i)$ as
$\Rstar$-modules, and as a
minimal free resolution of $A$ over $\Rstar$ of the form \eqref{eq9}, 
we have
\begin{equation}
0\lra\Ext{2}{R}{R/J}{R}(-\socd{A}-3)
\xra{\ \Phg\ }R/J\xra{\ \pg\ } A\lra0,
\lb{eq46}
\end{equation}
where $\pg$ is the natural surjection from $R/J$ to $A=R/I$ and
$\Phg$ is the map induced from the embedding $I/J\hra R/J$. 
Moreover 
$\Xsub{1}\Xsub{2}=\Xsub{2}\Xsub{1}$ and
$\Ysub{1}\Ysub{2}=\Ysub{2}\Ysub{1}$
with the notation of Section \ref{gor}.
\end{lem}

\begin{pf}
We have an exact sequence
\begin{equation*}
0\lra I/J\xra{\ \Phg\ } R/J\xra{\ \pg\ }A\lra0
\end{equation*}
over $R$ with the natural surjection $\pg$ and the embedding $\Phg$.
Since $R/J$ is \CM\ with $0:_{R/J} x_3=0$, it is free over $\Rstar$.
Hence $I/J$ is also free over $\Rstar$. Moreover, 
since $\Jbar=(\gbar_1,\gbar_2,\gbar_3)=\Ibar$, the homomorphism
$\pg\ot\Rbarstar:R/J\ot\Rbarstar\lra\Abar$ is an isomorphism.
The image of $\Phg$ is therefore contained in 
$x_3(R/J)$ and the above free resolution is minimal over $\Rstar$.
Since $R/J$ is a module over $R$, 
the matrix $\Xsub{j}$ can be chosen so that
it represents the linear map $\tm x_j:(R/J)(-1)\lra R/J$ over $\Rstar$.
Hence the relations
$\Xsub{1}\Xsub{2}=\Xsub{2}\Xsub{1}$ and
$\Ysub{1}\Ysub{2}=\Ysub{2}\Ysub{1}$ hold.
Taking the dual over $\Rstar$ of the above free resolution, 
we obtain the exact sequence
\begin{equation}
0\lra (R/J)^\vee\xra{\ \tpose{\Phg}\ } (I/J)^\vee\lra 
\Ext{1}{\Rstar}{A}{\Rstar}\cong A(\socd{A}+1)\lra0
\lb{eq47}
\end{equation}
over $R$, where the structures of $(R/J)^\vee,\ (I/J)^\vee$ and
$\Ext{1}{\Rstar}{A}{\Rstar}$ over $R$ are induced from those of
$R/J,\ I/J$ and $A$ respectively in a standard way.
Since $A$ is generated by a single element over $R$, so is $(I/J)^\vee$
by the minimal free resolution \eqref{eq47} over $\Rstar$. On the other hand, 
$\ann{(I/J)^\vee}\sset\ann{(R/J)^\vee}$,
since $\tpose{\Phg}$ is injective.
Besides, 
$J\sset\ann{(I/J)^\vee}$ and $\ann{(R/J)^\vee}
\sset\ann{(R/J)^\vee{}^\vee}=\ann{R/J}=J$. Hence $\ann{(I/J)^\vee}=J$,
that is, $(I/J)^\vee\cong (R/J)(\socd{A}+1)$.
In consequence, 
\begin{equation*}
I/J=(I/J)^\vee{}^\vee\cong\Ext{2}{R}{R/J}{R}(-\socd{A}-3),
\end{equation*}
which proves our assertion.
\end{pf}

\begin{thm}\lb{proc61}
Let $I$ be a homogeneous ideal in $R$ defining \Gor\ graded ring 
$A:=R/I$ of dimension zero. Suppose that $\ch{k}=0$. 
If $r=3$ and $\shmin{0:_Ax_3}\leq2$, then $A$ has the WLP.
\end{thm}

\begin{pf}
Assume that $r=3$ and that $\Gg$ is chosen sufficiently generally. 
We have $\shmin{\Ibar}\geq2$ since $\hei{\Ibar}=2$.
Besides, the condition $\shmin{0:_Ax_3}\leq2$ is equivalent to
$\shmin{\Ibar}\leq3$.
Let $(a;n^2_1\ddd n^2_a;n^3_1\ddd n^3_b)$ be the basic sequence of $I$.
By \Lem{\ref{proc5}}, it is enough to show that $n^3_1\geq n^2_a$.

\case{(i)}
Consider first the case $\shmin{\Ibar}=2$. With the use of 
\Cor{\ref{proc45}}, we see $\shmin{I}=3$. 
As we have already seen in \eqref{cond27} of \Lem{\ref{proc57}},
$n^3_1\geq n^2_a$. 

\par
Now assume that $\shmin{\Ibar}=3$ and let $g_i$, $\gbar_i$, $d_i$ and
$\gbar_{ij}$ be as stated just before \Lem{\ref{proc36}}.
We will show that we are led to a contradiction if $n^3_1 < n^2_a$. 

As in Section \ref{gor},
put $L:=\dss{i=1}{b}\Rstar(-n^3_b+n^3_i)$,
$F:=\dss{i=1}{b}\Rstar(-n^3_i)$, 
$\Lbar:=L/x_3L=\dss{i=1}{b}k(-n^3_b+n^3_i)$, and
$\Fbar:=F/x_3F=\dss{i=1}{b}k(-n^3_i)$, and
denote the canonical basis of $\Lbar$ (resp. $\Fbar$)
by $\{\ebar{}^L_1\ddd \ebar{}^L_{b}\}$ 
(resp. $\{\ebar{}^F_1\ddd \ebar{}^F_{b}\}$).
Suppose that $n^3_1<n^2_a$.
Put $p:=\max\{\ i\ |\ n^3_1-n^3_{b}+n^3_i\leq0,\ \ee{1}{i}{b}\ \}$
and $\pcheck:=\max\{\ i\ |\ n^3_i-n^3_{b}+n^3_p\leq0,\ \ee{1}{i}{b}\ \}$.
Since $n^3_b-n^3_1=n^2_a-1$ by \Cor{\ref{proc39}} and $n^3_1<n^2_a$ by
our assumption, we have $p\geq1$ and $\pcheck\geq1$.
Besides $p<b$, since $n^3_1>0$. 
Moreover, the sequence $-n^3_{b}+n^3_i$
$(\ee{1}{i}{b})$ ranges over all integers from $-n^3_{b}+n^3_1$ to $0$
as stated in \Lem{\ref{proc40}}, therefore $n^3_1=n^3_b-n^3_p$.
Put $q:=\max\{\ j\ |\ n^3_{p+1}=n^3_{p+j},\ \ee{1}{j}{b-p}\ \}$.
For every $\ee{p+1}{l}{p+q}$, we have $n^3_p < n^3_l$,
so that $n^3_b-n^3_p \geq n^3_b-n^3_l+1$. With the use of 
\Lem{\ref{proc58}}, we find that
\begin{align*}
&\gradn{\Ima{\Phg}}{n^3_b-n^3_{p+1}+2}\\
&\hspace{2em}\sset\gradn{\Ima{\Phg}}{\leq n^3_b-n^3_p+1}
\sset\Rstar[\Xsub{1},\Xsub{2}]
\gradn{\Ima{\Phg}}{\leq n^3_b-n^3_p}\quad\text{or}\\
&\thgsub{12l}\in
\Rstar[\Xsub{1},\Xsub{2}]
\gradn{\Ima{\Phg}}{\leq n^3_b-n^3_p}\quad
\text{for all $\ee{p+1}{l}{p+q}$}.
\end{align*}
Thus $\Rstar[\Xsub{1},\Xsub{2}]
\gradn{\Ima{\Phg}}{\leq n^3_1}
=\Rstar[\Xsub{1},\Xsub{2}]
\gradn{\Ima{\Phg}}{\leq n^3_b-n^3_p}\sset N_p$ by
\Cor{\ref{proc13}} or \Thm{\ref{proc18}}.
Observe that the polynomials 
$\gbar_1$ and $\gbar_2$ form an $\Rbar$ regular sequence
unless $\deg(\gcd(g_1,g_2))>0$, since $\Gg$ is sufficiently general.

\case{(ii)}
Suppose that $\gbar_1$ and $\gbar_2$ form an $\Rbar$-regular sequence or 
that $\gbar_{13}\neq0$. With the help of \eqref{eq15}, \eqref{eq22} and
\Lem{\ref{proc36}}, we find that
\begin{align*}
b&=\rk{k}{\Ext{2}{\Rbar}{\Abar}{\Rbar}}\\
&=\rk{k}{\Rbar(n^2_a+1)/\Imasup{\Rbar}{\tpose{\Psg}}\cap\Rbar(n^2_a+1)}
+\rk{k}{\Rbar(n)/(\gbar_{11},\gbar_{12},\gbar_{13})}\\
&<2\rk{k}{\Rbar(n^2_a+1)/\Imasup{\Rbar}{\tpose{\Psg}}\cap\Rbar(n^2_a+1)}.
\end{align*}
Let us consider $\Mbar_{3,1}$ defined in Section \ref{gor}. 
We have $\Mbar_{3,1}=
k[\Ybarsub{1},\Ybarsub{2}]\ebar^F_1
\sset\Fbar\cong
(0:_Ax_3)(-1)\cong\Ext{2}{\Rbar}{\Abar}{\Rbar}(-n^3_b-2)$.
Recall that the multiplication by $\Ybarsub{i}$ corresponds to
that by $x_i$ for $i=1,2$.
Let $e\neq0$ be a homogeneous element of 
$\Rbar(n^2_a+1)/\Imasup{\Rbar}{\tpose{\Psg}}\cap\Rbar(n^2_a+1)$ of
degree $-(n^2_a+1)$ which generate this module over $\Rbar$.
When considered in $\Ext{2}{\Rbar}{\Abar}{\Rbar}(-n^3_b-2)$, the
element $\ig(e)$ (see \eqref{eq22}) is of the minimal degree
since $-n>-(n^2_a+1)$ by \Lem{\ref{proc35}}, with 
$\deg(\ig(e))=-(n^2_a+1)+n^3_b+2=n^3_1$ 
by \Cor{\ref{proc39}}. Hence 
$\Mbar_{3,1}=\Rbar\cdot\ig(e)\cong
\Rbar(n^2_a+1)/\Imasup{\Rbar}{\tpose{\Psg}}\cap\Rbar(n^2_a+1)$.
On the other hand, $\Mbar_{3,1}\cong\Mbar_{2,1}$ by \eqref{eq16}, so that
\begin{equation*}
2\rk{k}{\Mbar_{3,1}}\leq \rk{k}{\Mbar_{3,\pcheck}}+
\rk{k}{\Mbar_{2,p}}\leq b
\end{equation*}
by \Lem{\ref{proc11}}, which is a contradiction. 

\case{(iii)}
Suppose that $\deg(\gcd(g_1,g_2))>0$ and that $\gbar_{13}=0$.
Let $J=(g_1,g_2,g_3)=(-h_2h_5,h_1h_5,h_2h_3-h_1h_4)$ be as in 
\Lem{\ref{proc59}}. Since $J$ has a minimal free resolution
\begin{equation*}
0\lra \dss{i=4}{5}R(-n^3_b-2+d_i)\xra{\ \Lg'\ }
\dss{i=1}{3}R(-d_i)\lra J\lra0
\end{equation*}
with 
$\Lg':=\begin{smallbmatrix}
h_3&h_1\\
h_4&h_2\\
h_5&0
\end{smallbmatrix}$
by \eqref{eq42} and since $\socd{A}=n^3_b-1$, there is an exact sequences
\begin{align*}
&0\lra (R/(h_5,h_2h_3-h_1h_4))(-d_4)\\
&\hspace{4em}\lra\Ext{2}{R}{R/J}{R}(-\socd{A}-3)
\lra (R/(h_1,h_2))(-d_5)\lra0.
\end{align*}
We consider $M_{1,1}$, $M_{2,1}$ and $M_{3,1}$ with the use of
\eqref{eq46}. In this case, $L=R/J$, $F=\Ext{2}{R}{R/J}{R}(-\socd{A}-3)$,
$\Xsub{1}\Xsub{2}=\Xsub{2}\Xsub{1}$ and
$\Ysub{1}\Ysub{2}=\Ysub{2}\Ysub{1}$.
The $\Rstar$-modules $L$, $L^\vee$ and $F$ are therefore modules 
over $R=\Rstar[x_1,x_2]$ through the surjections
$\Rstar[x_1,x_2]\thra\Rstar[\Xsub{1},\Xsub{2}]$,
$\Rstar[x_1,x_2]\thra\Rstar[\tpose{\Xsub{1}},\tpose{\Xsub{2}}]$
and $\Rstar[x_1,x_2]\thra\Rstar[\Ysub{1},\Ysub{2}]$ respectively.
So are $M_{1,1}$, $M_{2,1}$ and $M_{3,1}$.
Since $d_5=n^3_b-n+2 > n^3_b-(n^2_a+1)+2=n^3_1=d_4$
by \Cor{\ref{proc39}} and \Lems{\ref{proc35} and \ref{proc59}}, 
we find that $M_{3,1}\cong (R/(h_5,h_2h_3-h_1h_4))(-d_4)$ over $R$.
Since $L^\vee=(R/J)^\vee\cong\Ext{2}{R}{R/J}{R}(-2)=F(n^3_b)$ over $R$, 
we see $M_{2,1}\cong M_{3,1}(n^3_b)$. On the other hand,
since $M_{1,1}=\Phg M_{3,1}$ and $\Phg:F\lra L$
is injective, we have 
$(R/(h_1,h_2))(-d_5)\cong F/M_{3,1}\hra L/M_{1,1}$. 
Notice that $x_3^{\ng_0}(L/M_{1,1})\sset\Phg(F/M_{3,1})$
for some $\ng_0 > 0$ since $\det(\Phg)\neq0$.
Hence $(L/M_{1,1})^\vee\hra(F/M_{3,1})^\vee\cong(R/(h_1,h_2))^\vee(d_5)$.
Moreover, there is an injective map 
$\hg:M_{2,1}\hra(L/M_{1,1})^\vee$ by \Lem{\ref{proc10}} and \Thm{\ref{proc18}}.
Since $(R/(h_1,h_2))^\vee\cong\Ext{2}{R}{R/(h_1,h_2)}{R}
\cong R/(h_1,h_2)$ up to shift in grading,
we find in this manner that there is an injective homomorphism
$\hg':R/(h_5,h_2h_3-h_1h_4)\lra R/(h_1,h_2)$ over $R$, forgetting shift
of gradings. Since $h_1\hg'(1)=h_2\hg'(1)=0$, it must hold that
$h_1,h_2\in(h_5,h_2h_3-h_1h_4)$. But, then, $h_1,h_2\in(h_5)$.
Since $\gcd(h_1,h_2)\in k^\ast$, the degree of $h_5$ must be zero,
which contradicts our assumption that $\deg(\gcd(g_1,g_2))>0$.
\end{pf}

\section{The basic sequence of a complete intersection of three 
homogeneous polynomials}
\lb{ci}

Let $g_1,g_2,g_3\in R:=k[x_1,x_2,x_3]$ be homogeneous polynomials of degrees 
$d_1,d_2,d_3$ respectively with $2\leq d_1\leq d_2\leq d_3$ 
which form an $R$-regular sequence and let $I:=(g_1,g_2,g_3)\sset R$.
In this section, we give an explicit description of $\BR{I}$ for the
case $d_3 < d_1+d_2-1$.
When $d_3\geq d_1+d_2-1$, we have $\BR{I}\eqp B^0$ as have already been
proved in \eqref{cond27} of \Lem{\ref{proc57}}.

\begin{lem}\lb{proc30}
Let $B=(\ag;\ng^2_1\ddd\ng^2_\ag;\ng^3_1\ddd\ng^3_\bg)$ be a sequence
of integers such that $\ng^2_l=\ng^3_{l'}+1$ for some 
$\ee{1}{l}{\ag}$, $\ee{1}{l'}{\bg}$. Then
\begin{equation*}
B\sim (\ag;\ng^2_1\ddd\ng^2_{l-1},\ng^3_{l'},\ng^2_{l+1}\ddd\ng^2_{\ag};
\ng^3_1\ddd\ng^3_{l'-1},\ng^3_{l'+1}\ddd\ng^3_\bg).
\end{equation*}
\end{lem}

\begin{pf}
Our assertion follows from the equality
\begin{equation*}
\bnmm{\mg-\ng^2_l+1}{1}+\bnmm{\mg-\ng^3_{l'}}{0}
=\bnmm{\mg-\ng^3_{l'}}{1}+\bnmm{\mg-\ng^3_{l'}}{0}
=\bnmm{\mg-\ng^3_{l'}+1}{1}
\end{equation*}
for all $\mg\in\Zbf$.
\end{pf}

\par
We will denote the sequence $(l,l\ddd l)$ ($\mg$ times) by $l^\mg$ for 
$l\in\Zbf$, $\mg\in\Nbf$. 
When $d_3 < d_1+d_2-1$ and $d_1+d_2-d_3\geq 2c$ with $c\in\Nbf$,
we put
\begin{align*}
\ngbar^{c2}:=(&d_2,d_2+1\ddd d_3-1,(d_3)^2,(d_3+1)^2\ddd (d_3+c-1)^2,\\
&d_3+c,d_3+c+1\ddd d_2+d_1-c-1),\\
\ngbar^{c3}:=(&(d_2+d_1-1,d_2+d_1\ddd d_3+d_1-1),\\
&(d_2+d_1-2,d_2+d_1-1\ddd d_3+d_1),\\
&(d_2+d_1-3,d_2+d_1-2\ddd d_3+d_1+1),\\
&\cdots,\\
&(d_2+d_1-c,d_2+d_1-c+1\ddd d_3+d_1+c-2),\\
&((d_3+j,d_3+1+j\ddd d_3+d_1-1+j))_{\ee{c}{j}{d_2-1}}),
\end{align*}
and $B^c:=(d_1;\ngbar^{c2};\ngbar^{c3})$, where
\begin{align*}
\ngbar^{c2}:=(&(d_3)^2,(d_3+1)^2\ddd (d_3+c-1)^2,\\
&d_3+c,d_3+c+1\ddd d_3+d_1-c-1)
\end{align*}
if $d_2=d_3$. Note that $c\leq d_2-1$.

\begin{lem}\lb{proc31}
Suppose that $d_3 < d_1+d_2-1$ and that 
$d_1+d_2-d_3\geq 2c$ with $c\in\Nbf$.
Then $B^0 \sim B^c$.
\end{lem}

\begin{pf}
Since $d_3 < d_1+d_2-1$, the sequence $(d_3+1,d_3+2\ddd d_2+d_1-1)$
(resp. $(d_3,d_3+1\ddd d_2+d_1-2)$) is a subsequence of $\ngbar^{02}$
(resp. $\ngbar^{03}$).
Using \Lem{\ref{proc30}} $d_2+d_1-d_3-1$ times,
we find that $B^0\sim B^1$. If $c\geq2$, then $d_3+1 < d_1+d_2-2$,
so that we can repeat the same procedure as above to get $B^1\sim B^2$.
Hence $B^0\sim B^2$.
Proceeding in this way, we obtain $B^0\sim B^c$.
\end{pf}

\begin{thm}\lb{proc16}
Suppose that $d_3 < d_1+d_2-1$ and let $c$ be the unique
positive integer such that $2c\leq d_1+d_2-d_3\leq 2c+1$.
Then $\BR{I}\eqp B^c$.
\end{thm}

\begin{pf}
The maximal term of $\ngbar^{c2}$ is $d_2+d_1-c-1$ and the
minimal term of $\ngbar^{c3}$ is $d_3+c$. Their difference
$(d_3+c)-(d_2+d_1-c-1)$ is not negative by the choice of $c$, so that 
$d_2+d_1-c-1\leq d_3+c$. As for $\BR{I}$, one has $n^2_a\leq n^3_1$ 
by \Lem{\ref{proc5}} and \Thm{\ref{proc61}}. On the other hand $B^0\sim B^c$ by
\Lem{\ref{proc31}} and $\BR{I}\sim B^0$ by \Lem{\ref{proc32}}. 
Hence $\BR{I}\sim B^c$. Comparing $h_{\BR{I}}(\mg)$ and $h_{B^c}(\mg)$
$(\mg\in\Zbf)$, one finds that $\BR{I}\eqp B^c$.
\end{pf}

\end{document}